\RequirePackage{ifpdf}
\ifpdf 
\documentclass[pdftex]{sigma}
\else
\documentclass{sigma}
\fi

\newcommand{\intprod}{\;\rule{5pt}{.3pt}\rule{.3pt}{7pt}\;}
\newcommand{\End}{\operatorname{End}}
\newcommand{\Rho}{{\mathrm{P}}}
\def\verylongrightarrow{\relbar\joinrel\relbar\joinrel\longrightarrow}

\begin{document}


\renewcommand{\thefootnote}{$\star$}

\renewcommand{\PaperNumber}{060}

\FirstPageHeading

\ShortArticleName{The BGG Complex on Projective Space}

\ArticleName{The BGG Complex on Projective Space\footnote{This paper is a
contribution to the Special Issue ``Symmetry, Separation, Super-integrability and Special Functions~(S$^4$)''. The
full collection is available at
\href{http://www.emis.de/journals/SIGMA/S4.html}{http://www.emis.de/journals/SIGMA/S4.html}}}

\Author{Michael G.~EASTWOOD~$^\dag$ and A.~Rod GOVER~$^{\dag\ddag}$}

\AuthorNameForHeading{M.G.~Eastwood and A.R.~Gover}

\Address{$^\dag$~Mathematical Sciences Institute,
Australian National University, ACT 0200, Australia}
\EmailD{\href{mailto:meastwoo@member.ams.org}{meastwoo@member.ams.org}}

\Address{$^\ddag$~Department of Mathematics,
University of Auckland,\\
\phantom{$^\ddag$}~Private Bag 92019, Auckland 1142, New Zealand}
\EmailD{\href{mailto:r.gover@auckland.ac.nz}{r.gover@auckland.ac.nz}}


\ArticleDates{Received January 30, 2011, in f\/inal form June 18, 2011;  Published online June 23, 2011}

\Abstract{We give a complete construction of the Bernstein--Gelfand--Gelfand
complex on real or complex projective space using minimal ingredients.}

\Keywords{dif\/ferential complex; BGG complex; projective space; Lie algebra cohomology; parabolic geometry}

\Classification{53A20; 17B56; 58D19; 58J70}

\begin{flushright}
\begin{minipage}{80mm}
\it Dedicated to Willard Miller on his retirement from the University of
Minnesota.
\end{minipage}
\end{flushright}

\section{Introduction}
This article is concerned with dif\/ferential calculus on real or complex
projective space, in\-variant under projective transformations. These
transformations constitute a semisimple Lie group and projective space is hence
a homogeneous space of the form~$G/P$ for~$G$ semisimple. The subgroup~$P$ is
parabolic and, more generally, dif\/ferential geometries modelled on homogeneous
spaces of this form are known as `parabolic'~\cite{thebook}. Projective space
gives rise to projective dif\/ferential geometry in this sense. Conformal and CR
geometry are included amongst other examples of parabolic dif\/ferential
geometry. The interplay between the symmetries of projective space and its
invariant dif\/ferential operators is mediated by representation theory. In this
article we present the `BGG complex' on projective space as perhaps the
simplest of these constructions in one of the simplest of settings. We
anticipate that our approach will extend to $G/P$ in general and perhaps to
other homogeneous spaces. In~\cite{o}, Olver introduced complexes of
dif\/ferential operators on Euclidean space acting between `hyperforms'. These
are BGG complexes constructed directly, employing only af\/f\/ine invariance and
the associated Schur functors. He constructs some `easy' examples and observes
that ``Other examples, of greater complexity, can of course be constructed at
will, but the expressions rapidly get out of hand, even in low dimensional
spaces''. We maintain that the extra symmetry that the operators exhibit under
projective transformations allows one to control these expressions more
ef\/fectively. The fascinating combination of symmetries and dif\/ferential
operators is a def\/ining feature in the work of Willard Miller to whom we
dedicate this article.

The Bernstein--Gelfand--Gelfand (BGG) complexes on ${\mathbb{RP}}_n$ or
${\mathbb{CP}}_n$ are, by now, well-known complexes of vector bundles and
dif\/ferential operators between them generalising the de~Rham complex on
${\mathbb{RP}}_n$ and the holomorphic de~Rham complex on~${\mathbb{CP}}_n$,
respectively. An introduction to such complexes is given in~\cite{eVariations}
and the particular case of projective space is discussed in~\cite{eProjective}.
Usually, their construction involves choosing so-called `splitting
operators'~\cite{cd,css4} constructed from Kostant's Laplacian \cite{k} or
`quabla operator' \cite{cd} or from the Jantzen--Zuckerman translation
functor~\cite{v}. Here, we avoid the direct use of splitting operators, instead
relying only on diagram chasing, as is already done
in~\cite{eVariations,eProjective} in simple cases. In fact, the BGG complex of
holomorphic dif\/ferential operators on ${\mathbb{CP}}_2$ was already constructed
in this manner~\cite[p.~351]{eDuality} before it was realised by John
Rice~\cite{r} that complexes like this were dual to Lepowsky's construction
\cite{l} on the level of Verma modules. (On ${\mathbb{CP}}_1$ there is just one
family of BGG operators, already singled out for their invariance
in~\cite[Proposition~2.1]{et}.) In this article, for simplicity and cleanliness
we employ a spectral sequence to ef\/fect the diagram chasing. We employ
projective invariance to derive explicit formul{\ae} for the operators in the
projective BGG complex on ${\mathbb{RP}}_n$ in terms of the usual round metric
on the sphere.

As is often done in dif\/ferential geometry, when it is necessary to write out
tensors and their natural operations and we shall adorn them with upper or
lower indices corresponding to the tangent or cotangent bundle respectively.
For example, a vector f\/ield can be written as~$X^a$, a~one-form as $\omega_a$,
and the natural pairing between them as $X^a\omega_a$ in accordance with the
`Einstein summation convention'. For any tensor $\phi_{abc}$ we shall write its
symmetric part as $\phi_{(abc)}$ and its skew part as $\phi_{[abc]}$. For
example, to say that $\omega_{ab}$ is a two-form is to say that
\[\omega_{ab}=-\omega_{ba}\qquad
\mbox{or, equivalently,}\qquad\omega_{ab}=\omega_{[ab]}\qquad
\mbox{or, equivalently,}\qquad\omega_{(ab)}=0\]
and then
\[\nabla_{[a}\omega_{bc]}\qquad\mbox{and}\qquad
X^a\nabla_a\omega_{bc}-2(\nabla_{[b}X^a)\omega_{c]a},\]
for any torsion-free connection~$\nabla_a$, are the exterior derivative of
$\omega_{ab}$ and the Lie derivative of $\omega_{ab}$ in the direction of the
vector f\/ield~$X^a$, respectively. Such formul{\ae} are not meant to imply any
choice of local co\"ordinates. More precisely, this is Penrose's `abstract
index notation'~\cite{OT} and it formalises the conventions used by many
classical authors~-- see, for example, the discussion of projective dif\/ferential
geometry by Schouten~\cite{s}.

We shall view ${\mathbb{RP}}_n$ as a homogeneous space
\[{\mathbb{RP}}_n={\mathrm{SL}}(n+1,{\mathbb{R}})/P=G/P,\qquad\mbox{where}\quad
P=\left\{
\left\lgroup\begin{tabular}{c|ccc}$\ast$&$\ast$&$\!\cdots\!$&$\ast$\\ \hline
$0$\\
\raisebox{2pt}[15pt]{$\vdots$}&&
\raisebox{-2pt}[0pt][0pt]{\makebox[0pt]{\Huge$\ast$}}\\
$0$\end{tabular}\right\rgroup\right\}\]
or as a quotient of the `projective sphere'
\[S^n={\mathrm{SL}}(n+1,{\mathbb{R}})/P=G/P,\qquad\mbox{where}\quad
P=\left\{
\left\lgroup\begin{tabular}{c|ccc}$\lambda$&$\ast$&$\!\cdots\!$&$\ast$\\ \hline
$0$\\
\raisebox{2pt}[15pt]{$\vdots$}&&
\raisebox{-2pt}[0pt][0pt]{\makebox[0pt]{\Huge$\ast$}}\\
$0$\end{tabular}\right\rgroup\mbox{ s.t.\ }\lambda>0\right\}\]
under the antipodal map. In either case, it is convenient to write the
associated Lie algebra as
\begin{gather}\label{graded}
{\mathfrak{sl}}(n+1,{\mathbb{R}})={\mathfrak{g}}=
{\mathfrak{g}}_-\oplus{\mathfrak{g}}_0\oplus{\mathfrak{g}}_+,
\end{gather}
where{\samepage
\begin{gather*}
{\mathfrak{g}}_-=\left\{\!
\left\lgroup\!\begin{tabular}{c|ccc}$0$&$0$&$\!\cdots\!$&$0$\\ \hline
$\ast$\\
\raisebox{2pt}[15pt]{$\vdots$}&&
\raisebox{-2pt}[0pt][0pt]{\makebox[0pt]{\Huge$0$}}\\
$\ast$\end{tabular}\!\right\rgroup\!\right\},\!\quad
{\mathfrak{g}}_0=\left\{\!
\left\lgroup\!\begin{tabular}{c|ccc}$\ast$&$0$&$\!\cdots\!$&$0$\\ \hline
$0$\\
\raisebox{2pt}[15pt]{$\vdots$}&&
\raisebox{-2pt}[0pt][0pt]{\makebox[0pt]{\Huge$\ast$}}\\
$0$\end{tabular}\!\right\rgroup\!\right\},\!\quad
{\mathfrak{g}}_+=\left\{\!
\left\lgroup\!\begin{tabular}{c|ccc}$0$&$\ast$&$\!\cdots\!$&$\ast$\\ \hline
$0$\\
\raisebox{2pt}[15pt]{$\vdots$}&&
\raisebox{-2pt}[0pt][0pt]{\makebox[0pt]{\Huge$0$}}\\
$0$\end{tabular}\!\right\rgroup\!\right\}
\end{gather*}
and then ${\mathfrak{p}}={\mathfrak{g}}_0\oplus{\mathfrak{g}}_+$.}

From now on we shall discuss only real projective space ${\mathbb{RP}}_n$ or
its double cover, the sphere~$S^n$. The complex case is completely parallel
with real numbers being replaced by complex numbers everywhere and by working
in the holomorphic category rather than the smooth.

\section{An outline of the construction}\label{outline}

If ${\mathbb{V}}$ is a f\/inite-dimensional representation of~$P$, we shall
denote by $V$ the induced homogeneous vector bundle on~$G/P$ constructed as
\[V=G\times_P{\mathbb{V}}=G\times{\mathbb{V}}/\sim,
\qquad\mbox{where }(g,v)\sim\big(g p,p^{-1}v\big),\quad \forall\, p\in P.\]
Notice that if ${\mathbb{V}}$ is actually a $G$-module restricted to~$P$, then
$V$ is canonically trivialised
\begin{equation}\label{twist}
V=G\times_P{\mathbb{V}}\cong G/P\times{\mathbb{V}}\qquad\mbox{by}\quad
(g,v)\mapsto(gP,g v)\end{equation}
as a vector bundle (but not as a homogeneous vector bundle). Hence, in this
case $V$ is naturally equipped with a $G$-equivariant f\/lat connection~$\nabla$
obtained by transporting to $V$ the exterior derivative
$d:\Lambda^0\otimes{\mathbb{V}}\to\Lambda^1\otimes{\mathbb{V}}$ with
values in~${\mathbb{V}}$. More generally, the coupled de~Rham sequence
\begin{equation}\label{coupled}V\xrightarrow{\,\nabla\,}\Lambda^1\otimes V
\xrightarrow{\,\nabla\,}\Lambda^2\otimes V
\xrightarrow{\,\nabla\,}\Lambda^3\otimes V \xrightarrow{\,\nabla\,}\cdots
\xrightarrow{\,\nabla\,}\Lambda^{n-1}\otimes V
\xrightarrow{\,\nabla\,}\Lambda^n\otimes V\to 0\end{equation} is exact on the
level of germs and provides a resolution of ${\mathbb{V}}$ as a locally
constant sheaf on~$G/P$.

This general reasoning holds on any homogeneous space $G/P$ but the following
discussion is specif\/ic to ${\mathbb{RP}}_n$ or~$S^n$. Suppose ${\mathbb{V}}$ is
irreducible as a $G$-module. In this case we shall see that as a $P$-module
${\mathbb{V}}$ is f\/iltered
\begin{equation}\label{filtration}
{\mathbb{V}}={\mathbb{V}}_0+{\mathbb{V}}_1+{\mathbb{V}}_2+
\cdots+{\mathbb{V}}_{N-1}+{\mathbb{V}}_N\end{equation} meaning that
these are the subquotients listed in a natural order, starting on the
left with the smallest quotient of $\mathbb{V}$.  (In other words
${\mathbb{V}}_N$ is the smallest $P$-submodule in the f\/iltration, the
quotient ${\mathbb{V}}/{\mathbb{V}}_N$ has a f\/iltration
${\mathbb{V}}_0+{\mathbb{V}}_1+{\mathbb{V}}_2+
\cdots+{\mathbb{V}}_{N-1}$, and the meaning is now clear by induction.)
It follows that the bundle $V$ is correspondingly f\/iltered
\begin{equation}\label{correspondingfiltration}
V=V_0+V_1+V_2+\cdots+V_{N-1}+V_N\end{equation}
and now we claim that the connection $\nabla:V\to\Lambda^1\otimes V$, and
consequently the whole complex~(\ref{coupled}), is compatible with this
f\/iltration (as detailed in Theorem~\ref{compatible} below). The spectral
sequence of a~f\/iltered complex~\cite{c} now comes into play, having as its
$E_0$-level the following.
\begin{equation}\label{E0}\raisebox{-70pt}{\begin{picture}(300,148)(0,10)
\put(0,20){\vector(1,0){100}}
\put(0,20){\vector(0,1){80}}
\put(93,14){\makebox(0,0){$p$}}
\put(-4,91){\makebox(0,0){$q$}}
\put(10,150){\makebox(0,0){$V_0$}}
\put(50,150){\makebox(0,0){$\Lambda^1\otimes V_0$}}
\put(100,150){\makebox(0,0){$\Lambda^2\otimes V_0$}}
\put(150,150){\makebox(0,0){$\Lambda^3\otimes V_0$}}
\put(200,150){\makebox(0,0){$\Lambda^4\otimes V_0$}}
\put(250,150){\makebox(0,0){$\Lambda^5\otimes V_0$}}
\put(290,150){\makebox(0,0){$\cdots$}}
\put(290,120){\makebox(0,0){$\cdots$}}
\put(290,90){\makebox(0,0){$\cdots$}}
\put(290,60){\makebox(0,0){$\cdots$}}
\put(290,30){\makebox(0,0){$\cdots$}}
\put(50,135){\makebox(0,0){$\uparrow\scriptstyle\partial$}}
\put(100,135){\makebox(0,0){$\uparrow\scriptstyle\partial$}}
\put(150,135){\makebox(0,0){$\uparrow\scriptstyle\partial$}}
\put(200,135){\makebox(0,0){$\uparrow\scriptstyle\partial$}}
\put(250,135){\makebox(0,0){$\uparrow\scriptstyle\partial$}}
\put(50,120){\makebox(0,0){$V_1$}}
\put(100,120){\makebox(0,0){$\Lambda^1\otimes V_1$}}
\put(150,120){\makebox(0,0){$\Lambda^2\otimes V_1$}}
\put(200,120){\makebox(0,0){$\Lambda^3\otimes V_1$}}
\put(250,120){\makebox(0,0){$\Lambda^4\otimes V_1$}}
\put(100,105){\makebox(0,0){$\uparrow\scriptstyle\partial$}}
\put(150,105){\makebox(0,0){$\uparrow\scriptstyle\partial$}}
\put(200,105){\makebox(0,0){$\uparrow\scriptstyle\partial$}}
\put(250,105){\makebox(0,0){$\uparrow\scriptstyle\partial$}}
\put(100,90){\makebox(0,0){$V_2$}}
\put(150,90){\makebox(0,0){$\Lambda^1\otimes V_2$}}
\put(200,90){\makebox(0,0){$\Lambda^2\otimes V_2$}}
\put(250,90){\makebox(0,0){$\Lambda^3\otimes V_2$}}
\put(150,75){\makebox(0,0){$\uparrow\scriptstyle\partial$}}
\put(200,75){\makebox(0,0){$\uparrow\scriptstyle\partial$}}
\put(250,75){\makebox(0,0){$\uparrow\scriptstyle\partial$}}
\put(150,60){\makebox(0,0){$V_3$}}
\put(200,60){\makebox(0,0){$\Lambda^1\otimes V_3$}}
\put(250,60){\makebox(0,0){$\Lambda^2\otimes V_3$}}
\put(200,45){\makebox(0,0){$\uparrow\scriptstyle\partial$}}
\put(250,45){\makebox(0,0){$\uparrow\scriptstyle\partial$}}
\put(200,30){\makebox(0,0){$V_4$}}
\put(250,30){\makebox(0,0){$\Lambda^1\otimes V_4$}}
\end{picture}}\end{equation}
Here, the precise positioning of the co\"ordinate axes is a matter of
convention. The important property of the $E_0$-level is that the
dif\/ferentials $\partial$ are simply homomorphisms of vector bundles and we
shall show that they are induced by a complex of $G_0$-modules
\[{\mathbb{V}}\stackrel{\partial}{\longrightarrow}
{\mathfrak{g}}_-^*\otimes{\mathbb{V}}\stackrel{\partial}{\longrightarrow}
\Lambda^2{\mathfrak{g}}_-^*\otimes{\mathbb{V}}
\stackrel{\partial}{\longrightarrow}
\Lambda^3{\mathfrak{g}}_-^*\otimes{\mathbb{V}}
\stackrel{\partial}{\longrightarrow}\cdots\stackrel{\partial}{\longrightarrow}
\Lambda^{n-1}{\mathfrak{g}}_-^*\otimes{\mathbb{V}}
\stackrel{\partial}{\longrightarrow}
\Lambda^n{\mathfrak{g}}_-^*\otimes{\mathbb{V}},\]
\newpage

\noindent
where
\[G_0=\left\{
\left\lgroup\begin{tabular}{c|ccc}$\ast$&$0$&$\!\cdots\!$&$0$\\ \hline
$0$\\
\raisebox{2pt}[15pt]{$\vdots$}&&
\raisebox{-2pt}[0pt][0pt]{\makebox[0pt]{\Huge$\ast$}}\\
$0$\end{tabular}\right\rgroup\in{\mathrm{SL}}(n+1,{\mathbb{R}})\right\}.\]
Furthermore, we shall show that this complex def\/ines the Lie algebra
cohomology $H^r({\mathfrak{g}}_-,{\mathbb{V}})$, which in turn has been
computed by Kostant~\cite{k}. It follows that the $E_1$-level of the spectral
sequence is rather sparse, typically
\begin{equation}\label{sparse}\raisebox{-70pt}{\begin{picture}(300,148)(0,10)
\put(0,20){\vector(1,0){100}}
\put(0,20){\vector(0,1){80}}
\put(93,14){\makebox(0,0){$p$}}
\put(-4,91){\makebox(0,0){$q$}}
\put(10,150){\makebox(0,0){$H^0$}}
\put(50,150){\makebox(0,0){$0$}}
\put(100,150){\makebox(0,0){$0$}}
\put(150,150){\makebox(0,0){$0$}}
\put(200,150){\makebox(0,0){$0$}}
\put(250,150){\makebox(0,0){$0$}}
\put(290,150){\makebox(0,0){$\cdots$}}
\put(290,120){\makebox(0,0){$\cdots$}}
\put(290,90){\makebox(0,0){$\cdots$}}
\put(290,60){\makebox(0,0){$\cdots$}}
\put(290,30){\makebox(0,0){$\cdots$}}
\put(50,120){\makebox(0,0){$0$}}
\put(100,120){\makebox(0,0){$H^1$}}
\put(150,120){\makebox(0,0){$H^2$}}
\put(200,120){\makebox(0,0){$0$}}
\put(250,120){\makebox(0,0){$0$}}
\put(100,90){\makebox(0,0){$0$}}
\put(150,90){\makebox(0,0){$0$}}
\put(200,90){\makebox(0,0){$0$}}
\put(250,90){\makebox(0,0){$H^3$}}
\put(150,60){\makebox(0,0){$0$}}
\put(200,60){\makebox(0,0){$0$}}
\put(250,60){\makebox(0,0){$0$}}
\put(200,30){\makebox(0,0){$0$}}
\put(250,30){\makebox(0,0){$0$}}
\end{picture}}\end{equation}
where $H^0=V_0$ and, in particular, there is precisely one irreducible bundle
in each diagonal~$E_1^{p,d-p}$ for~$d$ f\/ixed. Since~(\ref{coupled})
resolves~${\mathbb{V}}$, we know that this spectral sequence is also converging
to ${\mathbb{V}}$ and the only way that this can happen is if the dif\/ferentials
f\/it together as a~resolution
\[0\to{\mathbb{V}}
\to H^0\to H^1\to H^2\to H^3\to\cdots\to H^{n-1}\to H^n\to 0.\]
This is the required BGG resolution. The rest of the article is devoted to
f\/illing in the details of this argument.

\section[The filtering of ${\mathbb{V}}$ as a $P$-module]{The f\/iltering of $\boldsymbol{{\mathbb{V}}}$ as a $\boldsymbol{P}$-module}\label{filtering}

Recall the decomposition (\ref{graded}) of
${\mathfrak{g}}={\mathfrak{sl}}(n+1,{\mathbb{R}})$ and consider the element
\[H=\frac{1}{n+1}
\left\lgroup\begin{tabular}{c|ccc}$n$&$0$&$\!\cdots\!$&$0$\\ \hline
$0$&$-1$&&$0$\\
\raisebox{2pt}[15pt]{$\vdots$}&&
\raisebox{1pt}[0pt][0pt]{\makebox[0pt]{\Large$\ddots$}}\\
$0$&$0$&&$-1$\end{tabular}\right\rgroup\in{\mathfrak{g}}_0.\]
uniquely characterised as lying in the centre of ${\mathfrak{g}}_0$ with
$[H,X]=X$, for $X\in{\mathfrak{g}}_+$. It is called the {\em grading
element\/}~\cite{thebook} of the $|1|$-graded Lie algebra~(\ref{graded}). If
the $G$-module ${\mathbb{V}}$ is restricted to~$G_0$, then it splits into
eigenspaces under~$H$. For the standard representation by matrix multiplication
on column vectors for example,
\begin{equation}\label{Hstandard}
H\left\lgroup\begin{array}{c}0\\ v_1\\ \vdots\\ v_n\end{array}\right\rgroup=
-\frac{1}{n+1}
\left\lgroup\begin{array}{c}0\\ v_1\\ \vdots\\ v_n\end{array}\right\rgroup
\qquad\mbox{and}\qquad
H\left\lgroup\begin{array}{c}x\\ 0\\ \vdots\\ 0\end{array}\right\rgroup=
\frac{n}{n+1}
\left\lgroup\begin{array}{c}x\\ 0\\ \vdots\\ 0\end{array}\right\rgroup.
\end{equation}
Rather than use the actual eigenvalues, which are rational in general, let us
subtract the lowest eigenvalue and write
\begin{equation}\label{splits}
{\mathbb{V}}={\mathbb{V}}_0\oplus{\mathbb{V}}_1\oplus{\mathbb{V}}_2
\oplus\cdots\oplus{\mathbb{V}}_{N-1}\oplus{\mathbb{V}}_N\end{equation}
for the eigenspace decomposition, noting that ${\mathfrak{g}}_+$ acts by
${\mathbb{V}}_j\to{\mathbb{V}}_{j+1}$ for all~$j$. It follows that
\begin{equation}\label{submodules}
{\mathbb{V}}^j\equiv{\mathbb{V}}_j\oplus{\mathbb{V}}_{j+1}\oplus\cdots\oplus
{\mathbb{V}}_N\end{equation}
are $P$-submodules of ${\mathbb{V}}$ for all $j$ and we have our
f\/iltration~(\ref{filtration}).

\section[The filtered complex $\Lambda^\bullet\otimes V$ and its spectral sequence]{The f\/iltered complex $\boldsymbol{\Lambda^\bullet\otimes V}$ and its spectral sequence}

The f\/iltration (\ref{submodules}) of ${\mathbb{V}}$ as a $P$-module certainly
induces a f\/iltration
\[
V=V^0\supseteq V^1\supseteq V^2\supseteq\cdots\supseteq V^{N-1}\supset V^N
\]
of $V$ by $G$-homogeneous vector bundles on~$G/P$.
\begin{theorem}\label{compatible} The connection
$\nabla:V\to\Lambda^1\otimes V$ is compatible with this filtration in the
sense that there is a commutative diagram
\[\begin{array}{ccc}V&\stackrel{\nabla}{\longrightarrow}&\Lambda^1\otimes V\\
\rule[.7pt]{.4pt}{6pt}\hspace*{.7pt}\cup
&&\rule[.7pt]{.4pt}{6pt}\hspace*{.7pt}\cup\\
V^k&\stackrel{\nabla}{\longrightarrow}&\Lambda^1\otimes V^{k-1}
\end{array}\]
for all $k=1,2,\ldots,N$.
\end{theorem}
\begin{proof}
Because the assertion is local, without loss of generality it suf\/f\/ices to prove
it on the standard af\/f\/ine co\"ordinate patch, namely
\begin{gather}\label{Q}
{\mathbb{R}}^n=Q/G_0\subset G/P={\mathbb{RP}}_n,\qquad\mbox{where}\quad
Q=\left\{
\left\lgroup\begin{tabular}{c|ccc}$\ast$&$0$&$\!\cdots\!$&$0$\\ \hline
$\ast$\\
\raisebox{2pt}[15pt]{$\vdots$}&&
\raisebox{-2pt}[0pt][0pt]{\makebox[0pt]{\Huge$\ast$}}\\
$\ast$\end{tabular}\right\rgroup\in{\mathrm{SL}}(n+1,{\mathbb{R}})\right\}.\!\!\!\!
\end{gather}
Recall (\ref{splits}) that ${\mathbb{V}}$ splits as a $G_0$-module. Hence the
same is true of $V$ restricted to this patch:
\begin{equation}\label{Vsplitting}V|_{{\mathbb{R}}^n}
=V_0\oplus V_1\oplus V_2\oplus\cdots\oplus V_{N-1}\oplus V_N.\end{equation}
To proceed we need a formula for~$\nabla$. The appendix discusses various
natural constructions on a general Lie group $G$, which we now specialise to
be the Abelian Lie group
\[{\mathbb{R}}^n=
G_-=\left\{
\left\lgroup\begin{tabular}{c|ccc}$1$&$0$&$\!\cdots\!$&$0$\\ \hline
$\ast$\\
\raisebox{2pt}[15pt]{$\vdots$}&&
\raisebox{-2pt}[0pt][0pt]{\makebox[0pt]{\LARGE${\mathrm{Id}}$}}\\
$\ast$\end{tabular}\right\rgroup\in{\mathrm{SL}}(n+1,{\mathbb{R}})\right\}.\]
and according to (\ref{nabla}) we f\/ind that $\nabla=d+\theta$, where
\begin{itemize}\itemsep=0pt
\item $V|_{{\mathbb{R}}^n}$ is trivialised by
$\begin{array}{ccc}G_-\times{\mathbb{V}}&
\stackrel{\simeq\:\:}{\rightarrow}&Q\times_{G_0}{\mathbb{V}}\\
\cap&\raisebox{-2pt}{$\nearrow$}&\|\\
Q\times{\mathbb{V}}&&V|_{{\mathbb{R}}^n}
\end{array}$ to def\/ine
$d:V|_{{\mathbb{R}}^n}\to\Lambda^1\otimes V|_{{\mathbb{R}}^n}$,
\item $\theta:V|_{{\mathbb{R}}^n}\to\Lambda^1\otimes V|_{{\mathbb{R}}^n}$ is
def\/ined by the same trivialisation together with the Maurer--Cartan form
$\theta$ on $G_-$; more specif\/ically,
\[V|_{{\mathbb{R}}^n}\cong
G_-\times{\mathbb{V}}\xrightarrow{\,\theta\otimes{\mathrm{Id}}\,}
\Lambda^1\otimes{\mathfrak{g}}_-\otimes{\mathbb{V}}
\xrightarrow{\,{\mathrm{Id}}\otimes\rho\,}\Lambda^1\otimes{\mathbb{V}}\cong
\Lambda^1\otimes V|_{{\mathbb{R}}^n},\]
where $\rho:{\mathfrak{g}}_-\otimes{\mathbb{V}}\to{\mathbb{V}}$ is the
representation of ${\mathfrak{g}}$ on ${\mathbb{V}}$ restricted
to~${\mathfrak{g}}_-$.
\end{itemize}
Evidently, $d$ preserves the splitting (\ref{Vsplitting}) whilst $\theta$
sends $V_k$ to $V_{k-1}$ for all $k=1,2,\ldots, N$. In particular, $\nabla$
sends $V^k=V_k\oplus\cdots$ to
$\Lambda^1\otimes V^{k-1}=\Lambda^1\otimes V_{k-1}\oplus\cdots$, as required.
\end{proof}

\begin{corollary}\label{easypeasy}
The complex~\eqref{coupled} is compatible with the
filtration~\eqref{correspondingfiltration}, i.e.~$\nabla$ sends
$\Lambda^p\otimes V^k$ to $\Lambda^{p+1}\otimes V^{k-1}$.
\end{corollary}

\begin{proof} In fact, in the trivialisation $V|_{{\mathbb{R}}^n}\cong
G_-\times{\mathbb{V}}$ employed in the proof of Theorem~\ref{compatible}
\begin{gather*}
\Lambda^p\otimes V|_{{\mathbb{R}}^n}\cong\Lambda^p\otimes{\mathbb{V}}\ni
\omega\otimes v\stackrel{\nabla}{\longmapsto}
d\omega\otimes v+(-1)^p\omega\wedge(\theta\intprod\rho)v\in
\Lambda^{p+1}\otimes{\mathbb{V}}\cong\Lambda^{p+1}\otimes V|_{{\mathbb{R}}^n}
\end{gather*}
and the conclusion is manifest.\end{proof}

According to this corollary, we may now consider the spectral sequence of the
f\/iltered complex $\nabla:\Lambda^\bullet\otimes V$ on~$G/P$, the $E_0$-level
of which is
\begin{equation}\label{E0pq}
E_0^{p,q}=\Lambda^{p+q}\otimes V_{-q}\qquad\mbox{with dif\/ferential}\quad
\partial: \ E_0^{p,q}\to E_0^{p,q+1},\end{equation}
where we have chosen to normalise (\ref{E0}) by placing $V_0$ at the origin. By
construction, the bundle~$V_k$ on $G/P$ is the homogeneous bundle induced from
${\mathbb{V}}_k={\mathbb{V}}^k/{\mathbb{V}}^{k+1}$ as a $P$-module
(cf.~(\ref{submodules})). We already know that, as an eigenspace for the
grading element~$H$ in the centre of~${\mathfrak{g}}_0$, the vector space
${\mathbb{V}}_k$ is a~$G_0$-module. By regarding it as the quotient
${\mathbb{V}}^k/{\mathbb{V}}^{k+1}$ we are equivalently making~${\mathbb{V}}_k$
into a $P$-module by decreeing that $G_+$ act trivially. By construction, the
$E_0$-dif\/ferential is $G$-equivariant. Furthermore, its def\/inition
\[\begin{array}{@{}ccccccccc}
0&\to&\Lambda^{p+1}\otimes V^{k-2}&\to&\Lambda^{p+1}\otimes V^{k-1}&\to&
\Lambda^{p+1}\otimes V_{k-1}&\to&0\\
&&\nabla\uparrow\phantom{\nabla}&&\nabla\uparrow\phantom{\nabla}&&
\partial\uparrow\phantom{\partial}\\
0&\to&\Lambda^p\otimes V^{k-1}&\to&\Lambda^p\otimes V^k&\to&
\Lambda^p\otimes V_k&\to&0
\end{array}\]
and the Leibniz rule ensure that it is linear over the functions. In other
words $\partial$ is a $G$-equivariant homomorphism of homogeneous bundles and,
as such, must be induced by a homomorphism of $P$-modules
$\Lambda^p({\mathfrak{g}}/{\mathfrak{p}})^*\otimes{\mathbb{V}}_k\to
\Lambda^{p+1}({\mathfrak{g}}/{\mathfrak{p}})^*\otimes{\mathbb{V}}_{k-1}$, which
we shall also denote by~$\partial$. In fact, from the formula for $\nabla$
displayed in the proof of Corollary~\ref{easypeasy}, we see that $\partial$ is
induced by
\begin{equation}\label{E0differential}
\Lambda^p({\mathfrak{g}}/{\mathfrak{p}})^*\otimes{\mathbb{V}}=
\Lambda^p{\mathfrak{g}}_-^*\otimes{\mathbb{V}}
\ni v_{\beta\gamma\cdots\delta}\stackrel{\partial}{\longmapsto}
\rho_{[\alpha}v_{\beta\gamma\cdots\delta]}\in
\Lambda^{p+1}{\mathfrak{g}}_-^*\otimes{\mathbb{V}}=
\Lambda^{p+1}({\mathfrak{g}}/{\mathfrak{p}})^*\otimes{\mathbb{V}},
\end{equation}
where ${\mathbb{V}}$ is regarded as a $P$-module by restricting the $G$ action
to $G_0$ and decreeing that $G_+$ act trivially. Since ${\mathfrak{g}}_-$ is
Abelian, this formula agrees with (\ref{thisispartial}) for the Koszul
dif\/ferential (used in~\cite{ce} to def\/ine Lie algebra cohomology). As observed
in the appendix, this is a complex of $G_0$-modules. Alternatively, we could
come to the same conclusion by restricting attention to a standard af\/f\/ine
co\"ordinate chart $G_-\cong{\mathbb{R}}^n\hookrightarrow{\mathbb{RP}}_n$ and
noticing that the kernel of
$d:\Lambda^p\otimes{\mathbb{V}}\to\Lambda^{p+1}\otimes{\mathbb{V}}$ consists
precisely of the left-invariant ${\mathbb{V}}$-valued $p$-forms under the
action of $G_-$ on itself. Since $\nabla=d+\partial$, we may now invoke the
second realisation of $H^r({\mathfrak{g}}_-,{\mathbb{V}})$ from
Theorem~\ref{geometricrealisation} in the Appendix.

To proceed we need to be explicit concerning the irreducible representation
${\mathbb{V}}$ of ${\mathrm{SL}}(n+1,{\mathbb{R}})$. In fact, one usually deals
with complex representations of ${\mathrm{SL}}(n+1,{\mathbb{C}})$ (starting
with ${\mathfrak{sl}}(n+1,{\mathbb{C}})$) but here there is no real dif\/ference
and we shall adopt the notation from~\cite{beastwood} in denoting such
representations by attaching non-negative integers to its Dynkin diagram
\[{\mathbb{V}}=\enskip\raisebox{-5pt}{\begin{picture}(140,17)
\put(0,5){\makebox(0,0){$\bullet$}}
\put(20,5){\makebox(0,0){$\bullet$}}
\put(40,5){\makebox(0,0){$\bullet$}}
\put(60,5){\makebox(0,0){$\bullet$}}
\put(80,5){\makebox(0,0){$\bullet$}}
\put(120,5){\makebox(0,0){$\bullet$}}
\put(140,5){\makebox(0,0){$\bullet$}}
\put(0,5){\line(1,0){90}}
\put(100,5){\makebox(0,0){\scriptsize$\cdots$}}
\put(110,5){\line(1,0){30}}
\put(0,12){\makebox(0,0){\scriptsize$a_1$}}
\put(20,12){\makebox(0,0){\scriptsize$a_2$}}
\put(40,12){\makebox(0,0){\scriptsize$a_3$}}
\put(60,12){\makebox(0,0){\scriptsize$a_4$}}
\put(80,12){\makebox(0,0){\scriptsize$a_5$}}
\put(100,12){\makebox(0,0){\scriptsize$\cdots$}}
\put(120,12){\makebox(0,0){\scriptsize$a_{n-1}$}}
\put(140,12){\makebox(0,0){\scriptsize$a_n$}}
\end{picture}}\]
(meaning that $-[a_1,a_2,\dots,a_n]$ is the lowest weight of this
representation with respect to the standard basis of fundamental weights).
With this notation, here is the conclusion of Kostant's computation~\cite{k}
of Lie algebra cohomology. We describe the result as an
${\mathrm{SL}}(n,{\mathbb{R}})$-module (by attaching non-negative integers to
an $A$-series Dynkin diagram with one fewer nodes) and f\/ix the action of $G_0$
by specifying how the grading element $H\in{\mathfrak{g}}_0$ acts.
\begin{equation}\label{cohomology}
\makebox[0pt][l]{\hspace*{-20pt}$\begin{array}{rcl}
H^0({\mathfrak{g}}_-,{\mathbb{V}})&=&
\enskip\raisebox{-5pt}{\begin{picture}(120,17)
\put(0,5){\makebox(0,0){$\bullet$}}
\put(20,5){\makebox(0,0){$\bullet$}}
\put(40,5){\makebox(0,0){$\bullet$}}
\put(60,5){\makebox(0,0){$\bullet$}}
\put(100,5){\makebox(0,0){$\bullet$}}
\put(120,5){\makebox(0,0){$\bullet$}}
\put(0,5){\line(1,0){70}}
\put(80,5){\makebox(0,0){\scriptsize$\cdots$}}
\put(90,5){\line(1,0){30}}
\put(0,12){\makebox(0,0){\scriptsize$a_2$}}
\put(20,12){\makebox(0,0){\scriptsize$a_3$}}
\put(40,12){\makebox(0,0){\scriptsize$a_4$}}
\put(60,12){\makebox(0,0){\scriptsize$a_5$}}
\put(80,12){\makebox(0,0){\scriptsize$\cdots$}}
\put(100,12){\makebox(0,0){\scriptsize$a_{n-1}$}}
\put(120,12){\makebox(0,0){\scriptsize$a_n$}}
\end{picture}}\qquad
H\leadsto -c\\[5pt]
H^1({\mathfrak{g}}_-,{\mathbb{V}})&=&
\enskip\quad\raisebox{-5pt}{\begin{picture}(130,17)
\put(0,5){\makebox(0,0){$\bullet$}}
\put(30,5){\makebox(0,0){$\bullet$}}
\put(50,5){\makebox(0,0){$\bullet$}}
\put(70,5){\makebox(0,0){$\bullet$}}
\put(110,5){\makebox(0,0){$\bullet$}}
\put(130,5){\makebox(0,0){$\bullet$}}
\put(0,5){\line(1,0){80}}
\put(90,5){\makebox(0,0){\scriptsize$\cdots$}}
\put(100,5){\line(1,0){30}}
\put(0,12){\makebox(0,0){\scriptsize$a_1+a_2+1$}}
\put(30,12){\makebox(0,0){\scriptsize$a_3$}}
\put(50,12){\makebox(0,0){\scriptsize$a_4$}}
\put(70,12){\makebox(0,0){\scriptsize$a_5$}}
\put(90,12){\makebox(0,0){\scriptsize$\cdots$}}
\put(110,12){\makebox(0,0){\scriptsize$a_{n-1}$}}
\put(130,12){\makebox(0,0){\scriptsize$a_n$}}
\end{picture}}\qquad
H\leadsto -c+a_1+1\\[5pt]
H^2({\mathfrak{g}}_-,{\mathbb{V}})&=&
\enskip\raisebox{-5pt}{\begin{picture}(140,17)
\put(0,5){\makebox(0,0){$\bullet$}}
\put(30,5){\makebox(0,0){$\bullet$}}
\put(60,5){\makebox(0,0){$\bullet$}}
\put(80,5){\makebox(0,0){$\bullet$}}
\put(120,5){\makebox(0,0){$\bullet$}}
\put(140,5){\makebox(0,0){$\bullet$}}
\put(0,5){\line(1,0){90}}
\put(100,5){\makebox(0,0){\scriptsize$\cdots$}}
\put(110,5){\line(1,0){30}}
\put(0,12){\makebox(0,0){\scriptsize$a_1$}}
\put(30,12){\makebox(0,0){\scriptsize$a_2+a_3+1$}}
\put(60,12){\makebox(0,0){\scriptsize$a_4$}}
\put(80,12){\makebox(0,0){\scriptsize$a_5$}}
\put(100,12){\makebox(0,0){\scriptsize$\cdots$}}
\put(120,12){\makebox(0,0){\scriptsize$a_{n-1}$}}
\put(140,12){\makebox(0,0){\scriptsize$a_n$}}
\end{picture}}\qquad
H\leadsto -c+a_1+a_2+2\\[5pt]
H^3({\mathfrak{g}}_-,{\mathbb{V}})&=&
\enskip\raisebox{-5pt}{\begin{picture}(140,17)
\put(0,5){\makebox(0,0){$\bullet$}}
\put(20,5){\makebox(0,0){$\bullet$}}
\put(50,5){\makebox(0,0){$\bullet$}}
\put(80,5){\makebox(0,0){$\bullet$}}
\put(120,5){\makebox(0,0){$\bullet$}}
\put(140,5){\makebox(0,0){$\bullet$}}
\put(0,5){\line(1,0){90}}
\put(100,5){\makebox(0,0){\scriptsize$\cdots$}}
\put(110,5){\line(1,0){30}}
\put(0,12){\makebox(0,0){\scriptsize$a_1$}}
\put(20,12){\makebox(0,0){\scriptsize$a_2$}}
\put(50,12){\makebox(0,0){\scriptsize$a_3+a_4+1$}}
\put(80,12){\makebox(0,0){\scriptsize$a_5$}}
\put(100,12){\makebox(0,0){\scriptsize$\cdots$}}
\put(120,12){\makebox(0,0){\scriptsize$a_{n-1}$}}
\put(140,12){\makebox(0,0){\scriptsize$a_n$}}
\end{picture}}\qquad
H\leadsto -c+a_1+a_2+a_3+3\\
\vdots\vdots\hspace*{20pt}&\vdots&\hspace*{20pt}\vdots\vdots\\
H^{n-1}({\mathfrak{g}}_-,{\mathbb{V}})&=&
\enskip\raisebox{-5pt}{\begin{picture}(140,17)
\put(0,5){\makebox(0,0){$\bullet$}}
\put(20,5){\makebox(0,0){$\bullet$}}
\put(40,5){\makebox(0,0){$\bullet$}}
\put(60,5){\makebox(0,0){$\bullet$}}
\put(100,5){\makebox(0,0){$\bullet$}}
\put(140,5){\makebox(0,0){$\bullet$}}
\put(0,5){\line(1,0){70}}
\put(80,5){\makebox(0,0){\scriptsize$\cdots$}}
\put(90,5){\line(1,0){50}}
\put(0,12){\makebox(0,0){\scriptsize$a_1$}}
\put(20,12){\makebox(0,0){\scriptsize$a_2$}}
\put(40,12){\makebox(0,0){\scriptsize$a_3$}}
\put(60,12){\makebox(0,0){\scriptsize$a_4$}}
\put(80,12){\makebox(0,0){\scriptsize$\cdots$}}
\put(100,12){\makebox(0,0){\scriptsize$a_{n-2}$}}
\put(140,12){\makebox(0,0){\scriptsize$a_{n-1}+a_n+1$}}
\end{picture}}\qquad\quad
H\leadsto -c+a_1+\cdots+a_{n-1}+n-1\\[5pt]
H^n({\mathfrak{g}}_-,{\mathbb{V}})&=&
\enskip\raisebox{-5pt}{\begin{picture}(130,17)
\put(0,5){\makebox(0,0){$\bullet$}}
\put(20,5){\makebox(0,0){$\bullet$}}
\put(40,5){\makebox(0,0){$\bullet$}}
\put(60,5){\makebox(0,0){$\bullet$}}
\put(100,5){\makebox(0,0){$\bullet$}}
\put(130,5){\makebox(0,0){$\bullet$}}
\put(0,5){\line(1,0){70}}
\put(80,5){\makebox(0,0){\scriptsize$\cdots$}}
\put(90,5){\line(1,0){40}}
\put(0,12){\makebox(0,0){\scriptsize$a_1$}}
\put(20,12){\makebox(0,0){\scriptsize$a_2$}}
\put(40,12){\makebox(0,0){\scriptsize$a_3$}}
\put(60,12){\makebox(0,0){\scriptsize$a_4$}}
\put(80,12){\makebox(0,0){\scriptsize$\cdots$}}
\put(100,12){\makebox(0,0){\scriptsize$a_{n-2}$}}
\put(130,12){\makebox(0,0){\scriptsize$a_{n-1}$}}
\end{picture}}\qquad
H\leadsto -c+a_1+\cdots+a_{n-1}+a_n+n
\end{array}$}\end{equation}
where
\[c=\frac{na_1+(n-1)a_2+(n-2)a_3+\cdots+2a_{n-1}+a_n}{n+1}\]
(obtained by acting on $[a_1,a_2,\dots,a_n]$ with the f\/irst column of the
inverse Cartan matrix for ${\mathfrak{sl}}(n+1)$). Notice that each cohomology
is an irreducible representation of~$G_0$.

\begin{theorem}\label{E1level}
The $E_1$-level of the spectral sequence of the filtered complex
$\nabla:\Lambda^\bullet\otimes V$ consists of irreducible homogeneous vector
bundles on~${\mathbb{RP}}_n$ under the action of
$G={\mathrm{SL}}(n+1,{\mathbb{R}})$. Only the following terms are non-zero
\[\begin{array}{@{}rcl}
E_1^{0,0}&\leftrightsquigarrow&H^0({\mathfrak{g}}_-,{\mathbb{V}})\\
E_1^{a_1+1,-a_1}&\leftrightsquigarrow&H^1({\mathfrak{g}}_-,{\mathbb{V}})\\
E_1^{a_1+a_2+2,-a_1-a_2}&\leftrightsquigarrow&
H^2({\mathfrak{g}}_-,{\mathbb{V}})\\
E_1^{a_1+a_2+a_3+3,-a_1-a_2-a_3}&\leftrightsquigarrow&
H^3({\mathfrak{g}}_-,{\mathbb{V}})\\
\vdots\vdots\hspace*{20pt}&\vdots&\hspace*{20pt}\vdots\vdots\\
E_1^{N+n-1-a_n,-N+a_n}\enskip=\enskip
E_1^{a_1+a_2+a_3+\cdots+a_{n-1}+n-1,-a_1-a_2-a_3-\cdots-a_{n-1}}
&\leftrightsquigarrow&
H^{n-1}({\mathfrak{g}}_-,{\mathbb{V}})\\
E_1^{N+n,-N}\enskip=\enskip
E_1^{a_1+a_2+a_3+\cdots+a_{n-1}+a_n+n,-a_1-a_2-a_3-\cdots-a_{n-1}-a_n}
&\leftrightsquigarrow&
H^n({\mathfrak{g}}_-,{\mathbb{V}})
\end{array}\]
meaning that the bundle in question is induced by the $G_0$-module as listed
and extended trivially as a $P$-module.
\end{theorem}

\begin{proof} According to (\ref{E0pq}), we already know that
$\Lambda^r\otimes V$ is spread along the
$r^{\mathrm{th}}$ diagonal
\[\Lambda^r\otimes V=\Lambda^r\otimes V_0+\Lambda^r\otimes
V_1+\Lambda^r\otimes V_2+\cdots =
E_0^{r,0}+E_0^{r+1,-1}+E_0^{r+2,-2}+\cdots\]
of the $E_0$-level of the spectral sequence and that the $E_0$-dif\/ferential is
induced by the Koszul dif\/ferential~(\ref{E0differential}) def\/ining the Lie
algebra cohomology $H^r({\mathfrak{g}}_-,{\mathbb{V}})$. We see from
(\ref{cohomology}) that each of these $H^r({\mathfrak{g}}_-,{\mathbb{V}})$ is
irreducible and so it follows that the $r^{\mathrm{th}}$ diagonal of the
$E_1$-level consists of a single irreducible homogeneous vector bundle. To
complete the proof it suf\/f\/ices to locate the position of this bundle along this
particular diagonal and, to do this, the action of the grading element $H$
turns out to be suf\/f\/icient. More precisely,
\[\begin{array}{@{}rclrcl}
H\mbox{ acts by}&-c&\mbox{on }V_0&
\therefore\;H\mbox{ acts by}&-c+k&\mbox{on }V_k,\\
H\mbox{ acts by}&-1&\mbox{on }{\mathfrak{g}}_-&
\therefore\;H\mbox{ acts by}&p&\mbox{on }\Lambda^p{\mathfrak{g}}_-^*,\\
&&&\therefore\;H\mbox{ acts by}&-c+k+p&\mbox{on }
\Lambda^p{\mathfrak{g}}_-^*\otimes{\mathbb{V}}_k,\\
&&&\therefore\;H\mbox{ acts by}&-c+p&\mbox{on }
\Lambda^{p+q}{\mathfrak{g}}_-^*\otimes{\mathbb{V}}_{-q}\leftrightsquigarrow
E_0^{p,q}.
\end{array}\]
Therefore $H$ acts by $-c+p$ on the $G_0$-module inducing $E_1^{p,q}$ and so,
from the action of $H$ in table~(\ref{cohomology}), the bundle induced by
$H^r({\mathfrak{g}}_-,{\mathbb{V}})$ is located at
$E_1^{a_1+a_2+\cdots+a_r+r,-a_1-a_2-\cdots-a_r}$.

Finally, we are claiming that $N=a_1+a_2+\cdots+a_n$. To see this we note that,
since $H$ acts on $V_0$ by $-c$, it acts on $V_N$ by $N-c$. On the other hand,
in accordance with the action of the longest element of the Weyl group,
\[{\mathbb{V}}^*=\enskip\raisebox{-5pt}{\begin{picture}(140,17)
\put(0,5){\makebox(0,0){$\bullet$}}
\put(20,5){\makebox(0,0){$\bullet$}}
\put(60,5){\makebox(0,0){$\bullet$}}
\put(80,5){\makebox(0,0){$\bullet$}}
\put(100,5){\makebox(0,0){$\bullet$}}
\put(120,5){\makebox(0,0){$\bullet$}}
\put(140,5){\makebox(0,0){$\bullet$}}
\put(0,5){\line(1,0){30}}
\put(40,5){\makebox(0,0){\scriptsize$\cdots$}}
\put(50,5){\line(1,0){90}}
\put(0,12){\makebox(0,0){\scriptsize$a_n$}}
\put(20,12){\makebox(0,0){\scriptsize$a_{n-1}$}}
\put(40,12){\makebox(0,0){\scriptsize$\cdots$}}
\put(60,12){\makebox(0,0){\scriptsize$a_5$}}
\put(80,12){\makebox(0,0){\scriptsize$a_4$}}
\put(100,12){\makebox(0,0){\scriptsize$a_3$}}
\put(120,12){\makebox(0,0){\scriptsize$a_2$}}
\put(140,12){\makebox(0,0){\scriptsize$a_1$}}
\end{picture}}\enskip={\mathbb{V}}_N^*+{\mathbb{V}}_{N-1}^*+\cdots+
{\mathbb{V}}_2^*+{\mathbb{V}}_1^*+{\mathbb{V}}_0^*\]
and so $H$ acts on ${\mathbb{V}}_N^*$ by $-c^\prime$, where
\[c^\prime=\frac{na_n+(n-1)a_{n-1}+\cdots+3a_3+2a_2+a_1}{n+1}.\]
We conclude that $N=c+c^\prime=a_1+a_2+a_3+\cdots+a_{n-1}+a_n$, as required.
\end{proof}

As outlined in Section~\ref{outline}, the $E_1$-level of this spectral sequence is
rather sparse~(\ref{sparse}) and Theorem~\ref{E1level} says exactly how sparse.
In particular, since there is only one non-zero bundle on each diagonal of the
$E_1$-level, the general theory of spectral sequences provides a complex of
dif\/ferential operators
\[H^0\to H^1\to H^2\to H^3\to\cdots\to H^{n-1}\to H^n\to 0\]
whose cohomology on the level of sheaves coincides with that of
(\ref{coupled}), namely ${\mathbb{V}}$ for $\ker:H^0\to H^1$ and otherwise
zero. The bundle $H^r$ is induced on ${\mathbb{RP}}_n=G/P$ from the
$G_0$-module $H^r({\mathfrak{g}}_-,{\mathbb{V}})$ extended trivially as as
$P$-module. This is the BGG resolution (constructed on a general $G/P$ by
Lepowsky~\cite{l} on the level of generalised Verma modules). More explicitly,
in the discussion following Corollary~\ref{easypeasy}, in any af\/f\/ine
co\"ordinate patch we identif\/ied
$\nabla:\Lambda^p\otimes V\to\Lambda^{p+1}\otimes V$ as $d+\partial$ where
$V|_{{\mathbb{R}}^n}$ is trivialised as ${\mathbb{R}}^n\times{\mathbb{V}}$ and
\[E_0^{p,q}=\Lambda^{p+q}\otimes V_{-q}\cong
\Lambda^{p+q}\otimes{\mathbb{V}}_{-q}\stackrel{d}{\longrightarrow}
\Lambda^{p+q+1}\otimes{\mathbb{V}}_{-q}\cong\Lambda^{p+q+1}\otimes V_{-q}=
E_0^{p+1,q}\]
is the exterior derivative with values in~${\mathbb{V}}_{-q}$. Adding these
operators to the diagram (\ref{E0}) gives a~double complex as employed by
Baston~\cite[p.~120]{b}. We conclude that in any af\/f\/ine co\"ordinate patch, our
spectral sequence (of a f\/iltered complex) coincides with Baston's spectral
sequence (of a double complex). In particular, the operators $H^k\to H^{k+1}$
are obtained as zigzag compositions of the f\/irst order dif\/ferential operators
$d$ together with choices of algebraic splittings of the operators~$\partial$
(Baston and others use the algebraic adjoint $\partial^*$ introduced by
Kostant~\cite{k}). This conf\/irms that the resulting operators $H^k\to H^{k+1}$
are dif\/ferential. By using the spectral sequence of a f\/iltered complex as we
have done, it is manifest that the dif\/ferential operators $H^k\to H^{k+1}$ are
independent of any choice of splittings and that the whole construction and
resulting BGG complex is $G$-equivariant. It is also clear that the f\/irst steps
in this approach can be taken on an arbitrary homogeneous space~$G/P$.

\section{Formul{\ae} for the BGG operators}
Already, the af\/f\/ine invariance of the operators occurring in the BGG resolution
on ${\mathbb{RP}}_n$ f\/ixes their formul{\ae} with respect to the f\/lat
connection on ${\mathbb{R}}^n\subset{\mathbb{RP}}_n$ as follows (af\/f\/ine
invariance being ensured by noting that $Q$ as in (\ref{Q}) is acting on
${\mathbb{R}}^n=Q/G_0$ by af\/f\/ine transformations). Firstly, consider the symbol
of the BGG dif\/ferential operator $H^r\to H^{r+1}$. As a homomorphism of
homogeneous bundles $\bigodot^s\!\Lambda^1\otimes H^r\to H^{r+1}$ for
some~$s$, it is induced by a homomorphism of $G_0$-modules
\[\textstyle
\bigodot^s\!{\mathfrak{g}}_-^*\otimes H^r({\mathfrak{g}}_-,{\mathbb{V}})\to
H^{r+1}({\mathfrak{g}}_-,{\mathbb{V}})\]
for some $s$, which we can determine just from the action of the grading
element. Specif\/ically, we know that
\[\begin{array}{@{}rcl}
H\mbox{ acts by}&-c+a_1+a_2+\cdots+a_r+r&
\mbox{on }H^r({\mathfrak{g}}_-,{\mathbb{V}})\\
H\mbox{ acts by}&-1&
\mbox{on }{\mathfrak{g}}_-\\
\therefore\enskip H\mbox{ acts by}&s&
\mbox{on }\bigodot^s\!{\mathfrak{g}}_-^*.
\end{array}\]
It is immediate that $s=a_{r+1}+1$. Furthermore, from the
Littlewood--Richardson rules~(e.g.,~\cite{fh}), as an
${\mathfrak{sl}}(n,{\mathbb{R}})$-module
$H^{r+1}({\mathfrak{g}}_-,{\mathbb{V}})$ occurs with multiplicity one in the
following decomposition
\[\begin{array}{@{}lcc}
\hspace*{21.4pt}\bigodot^{a_{r+1}+1}\!{\mathfrak{g}}_-^*
\otimes H^r({\mathfrak{g}}_-,{\mathbb{V}})
&=&\cdots\oplus H^{r+1}({\mathfrak{g}}_-,{\mathbb{V}})\oplus\cdots\\
\hspace*{79.5pt}\|&&\|\\
\raisebox{-5pt}{\begin{picture}(70,17)
\put(0,5){\makebox(0,0){$\bullet$}}
\put(30,5){\makebox(0,0){$\bullet$}}
\put(70,5){\makebox(0,0){$\bullet$}}
\put(0,5){\line(1,0){40}}
\put(50,5){\makebox(0,0){\scriptsize$\cdots$}}
\put(60,5){\line(1,0){10}}
\put(0,12){\makebox(0,0){\scriptsize$a_{r+1}+1$}}
\put(30,12){\makebox(0,0){\scriptsize$0$}}
\put(50,12){\makebox(0,0){\scriptsize$\cdots$}}
\put(70,12){\makebox(0,0){\scriptsize$0$}}
\end{picture}}\enskip\otimes\enskip
\raisebox{-5pt}{\begin{picture}(140,17)
\put(0,5){\makebox(0,0){$\bullet$}}
\put(60,5){\makebox(0,0){$\bullet$}}
\put(100,5){\makebox(0,0){$\bullet$}}
\put(140,5){\makebox(0,0){$\bullet$}}
\put(0,5){\line(1,0){10}}
\put(20,5){\makebox(0,0){\scriptsize$\cdots$}}
\put(30,5){\line(1,0){80}}
\put(120,5){\makebox(0,0){\scriptsize$\cdots$}}
\put(130,5){\line(1,0){10}}
\put(0,12){\makebox(0,0){\scriptsize$a_1$}}
\put(20,12){\makebox(0,0){\scriptsize$\cdots$}}
\put(60,12){\makebox(0,0){\scriptsize$a_r+a_{r+1}+1$}}
\put(100,12){\makebox(0,0){\scriptsize$a_{r+2}$}}
\put(120,12){\makebox(0,0){\scriptsize$\cdots$}}
\put(140,12){\makebox(0,0){\scriptsize$a_n$}}
\end{picture}}&=&
\oplus\enskip\raisebox{-5pt}{\begin{picture}(140,17)
\put(0,5){\makebox(0,0){$\bullet$}}
\put(40,5){\makebox(0,0){$\bullet$}}
\put(80,5){\makebox(0,0){$\bullet$}}
\put(140,5){\makebox(0,0){$\bullet$}}
\put(0,5){\line(1,0){10}}
\put(20,5){\makebox(0,0){\scriptsize$\cdots$}}
\put(30,5){\line(1,0){80}}
\put(120,5){\makebox(0,0){\scriptsize$\cdots$}}
\put(130,5){\line(1,0){10}}
\put(0,12){\makebox(0,0){\scriptsize$a_1$}}
\put(20,12){\makebox(0,0){\scriptsize$\cdots$}}
\put(40,12){\makebox(0,0){\scriptsize$a_r$}}
\put(80,12){\makebox(0,0){\scriptsize$a_{r+1}+a_{r+2}+1$}}
\put(120,12){\makebox(0,0){\scriptsize$\cdots$}}
\put(140,12){\makebox(0,0){\scriptsize$a_n$}}
\end{picture}}\enskip\oplus
\end{array}\]
The symbol of the BGG operator $H^r\to H^{r+1}$
is thus determined uniquely (up to scale) as induced by the projection onto
this summand. Similarly, there are no lower order terms since no appropriate
invariant homomorphisms are available. We record this conclusion as follows.
\begin{theorem}
Each bundle $H^r$ in the BGG complex on ${\mathbb{RP}}_n$ is an irreducible
tensor bundle and the differential operator $\nabla^{(s)}:H^r\to H^{r+1}$ is
given by $\phi\mapsto\pi(\nabla^s\phi)$ where $\nabla$ is the flat affine
connection and $\pi:\bigodot^s\!\Lambda^1\otimes H^r\to H^{r+1}$ is the unique
projection onto this irreducible summand.
\end{theorem}
Equivalently, the bundles $H^r$ and the operators between them are exhibited as
Young tableau in~\cite{eProjective}. It is often useful, however, to be
able to write the BGG complex globally on ${\mathbb{RP}}_n$ without recourse to
af\/f\/ine co\"ordinates and projective invariance. For this, we shall use the
round metric on the sphere and the corresponding Levi Civita connection on
$S^n$ or~${\mathbb{RP}}_n$. By way of normalisation, if $g_{ab}$ denotes the
round metric and $\delta_a{}^b$ the identity endomorphism on the tangent
bundle, let us choose the radius of the sphere so that
\begin{equation}\label{cc}
R_{ab}{}^c{}_d=\delta_a{}^cg_{bd}-\delta_b{}^cg_{ad}\qquad\mbox{where}\quad
(\nabla_a\nabla_b-\nabla_b\nabla_a)X^d=R_{ab}{}^c{}_dX^d\end{equation}
def\/ines the Riemann curvature tensor. For any irreducible covariant tensor
bundle $E$ on ${\mathbb{RP}}_n$, let us write $\nabla^2$ for the composition
\[\textstyle\bigodot^r\!\Lambda^1\otimes E\xrightarrow{\,\nabla\circ\nabla\,}
\Lambda^1\otimes\Lambda^1\otimes\bigodot^r\!\Lambda^1\otimes E\to
\bigodot^{r+2}\!\Lambda^1\otimes E\]
and $g$ for the composition
\[\textstyle\bigodot^r\!\Lambda^1\otimes E
\xrightarrow{\,g\otimes{\mathrm{Id}}\,}
\bigodot^2\!\Lambda^1\otimes\bigodot^r\!\Lambda^1\otimes E\to
\bigodot^{r+2}\!\Lambda^1\otimes E.\]
\begin{theorem}\label{round_formulae}
If $s$ is odd, the operator $\nabla^{(s)}:H^r\to H^{r+1}$ is given by
\[\pi\big((\nabla^2+(s-1)^2g)\cdots(\nabla^2+16g)(\nabla^2+4g)\nabla\big)\]
and, if $s$ is even,
\[\pi\big((\nabla^2+(s-1)^2g)\cdots(\nabla^2+9g)(\nabla^2+g)\big).\]
\end{theorem}
\begin{proof}
It is convenient to use some notation from~\cite{eProjective} where the
projective BGG operators are written as acting between covariant tensor bundles
specif\/ied by weighted Young tableau: 
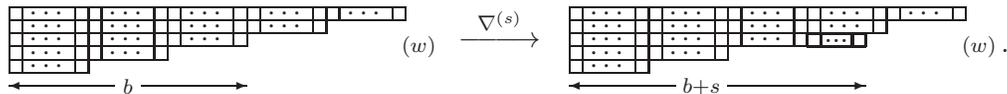
\[\begin{picture}(170,35)(0,85)
\put(0,115){\line(1,0){150}}
\put(0,110){\line(1,0){150}}
\put(0,105){\line(1,0){120}}
\put(0,100){\line(1,0){90}}
\put(0,95){\line(1,0){60}}
\put(0,90){\line(1,0){30}}
\put(0,90){\line(0,1){25}}
\put(30,90){\line(0,1){25}}
\put(60,95){\line(0,1){20}}
\put(90,100){\line(0,1){15}}
\put(120,105){\line(0,1){10}}
\put(150,110){\line(0,1){5}}
\put(15,92.5){\makebox(0,0){$\cdots$}}
\put(15,97.5){\makebox(0,0){$\cdots$}}
\put(15,102.5){\makebox(0,0){$\cdots$}}
\put(15,107.5){\makebox(0,0){$\cdots$}}
\put(15,112.5){\makebox(0,0){$\cdots$}}
\put(45,97.5){\makebox(0,0){$\cdots$}}
\put(45,102.5){\makebox(0,0){$\cdots$}}
\put(45,107.5){\makebox(0,0){$\cdots$}}
\put(45,112.5){\makebox(0,0){$\cdots$}}
\put(75,102.5){\makebox(0,0){$\cdots$}}
\put(75,107.5){\makebox(0,0){$\cdots$}}
\put(75,112.5){\makebox(0,0){$\cdots$}}
\put(105,107.5){\makebox(0,0){$\cdots$}}
\put(105,112.5){\makebox(0,0){$\cdots$}}
\put(135,112.5){\makebox(0,0){$\cdots$}}
\put(5,90){\line(0,1){25}}
\put(25,90){\line(0,1){25}}
\put(35,95){\line(0,1){20}}
\put(55,95){\line(0,1){20}}
\put(65,100){\line(0,1){15}}
\put(85,100){\line(0,1){15}}
\put(95,105){\line(0,1){10}}
\put(115,105){\line(0,1){10}}
\put(125,110){\line(0,1){5}}
\put(145,110){\line(0,1){5}}
\put(40,85){\vector(-1,0){40}}
\put(50,85){\vector(1,0){40}}
\put(45,85){\makebox(0,0){$\scriptstyle b$}}
\put(155,100){\makebox(0,0){$\scriptstyle (w)$}}
\end{picture}
\raisebox{15pt}{$\stackrel{\nabla^{(s)}}{\verylongrightarrow}\quad$}
\begin{picture}(170,35)(0,85)
\put(0,115){\line(1,0){150}}
\put(0,110){\line(1,0){150}}
\put(0,105){\line(1,0){120}}
\put(0,100){\line(1,0){90}}
\put(0,95){\line(1,0){60}}
\put(0,90){\line(1,0){30}}
\put(0,90){\line(0,1){25}}
\put(30,90){\line(0,1){25}}
\put(60,95){\line(0,1){20}}
\put(90,100){\line(0,1){15}}
\put(120,105){\line(0,1){10}}
\put(150,110){\line(0,1){5}}
\put(15,92.5){\makebox(0,0){$\cdots$}}
\put(15,97.5){\makebox(0,0){$\cdots$}}
\put(15,102.5){\makebox(0,0){$\cdots$}}
\put(15,107.5){\makebox(0,0){$\cdots$}}
\put(15,112.5){\makebox(0,0){$\cdots$}}
\put(45,97.5){\makebox(0,0){$\cdots$}}
\put(45,102.5){\makebox(0,0){$\cdots$}}
\put(45,107.5){\makebox(0,0){$\cdots$}}
\put(45,112.5){\makebox(0,0){$\cdots$}}
\put(75,102.5){\makebox(0,0){$\cdots$}}
\put(75,107.5){\makebox(0,0){$\cdots$}}
\put(75,112.5){\makebox(0,0){$\cdots$}}
\put(105,107.5){\makebox(0,0){$\cdots$}}
\put(105,112.5){\makebox(0,0){$\cdots$}}
\put(135,112.5){\makebox(0,0){$\cdots$}}
\put(5,90){\line(0,1){25}}
\put(25,90){\line(0,1){25}}
\put(35,95){\line(0,1){20}}
\put(55,95){\line(0,1){20}}
\put(65,100){\line(0,1){15}}
\put(85,100){\line(0,1){15}}
\put(95,105){\line(0,1){10}}
\put(115,105){\line(0,1){10}}
\put(125,110){\line(0,1){5}}
\put(145,110){\line(0,1){5}}
\put(40,85){\vector(-1,0){40}}
\put(60,85){\vector(1,0){52}}
\put(50,85){\makebox(0,0){$\scriptstyle b+s$}}
\put(155,100){\makebox(0,0){$\scriptstyle (w)\makebox[0pt][l]{$\,.$}$}}
\thicklines
\put(90,100){\line(1,0){22}}
\put(90,105){\line(1,0){22}}
\put(90,100){\line(0,1){5}}
\put(95,100){\line(0,1){5}}
\put(107,100){\line(0,1){5}}
\put(112,100){\line(0,1){5}}
\put(98,102.5){\makebox(0,0){.}}
\put(101,102.5){\makebox(0,0){.}}
\put(104,102.5){\makebox(0,0){.}}
\end{picture}\]
Here, if it is the $r^{\mathrm{th}}$ row to which the boxes on the right hand
side are being added and there are a total of $v$ boxes on the left hand side
(so that $v$ is the valence of the corresponding tensor), then $w+r=s+v+b$
(see~\cite{eProjective}). In
particular, if we consider only f\/irst order operators
\begin{equation}\label{firstorder}\begin{picture}(170,35)(0,85)
\put(0,115){\line(1,0){150}}
\put(0,110){\line(1,0){150}}
\put(0,105){\line(1,0){120}}
\put(0,100){\line(1,0){90}}
\put(0,95){\line(1,0){60}}
\put(0,90){\line(1,0){30}}
\put(0,90){\line(0,1){25}}
\put(30,90){\line(0,1){25}}
\put(60,95){\line(0,1){20}}
\put(90,100){\line(0,1){15}}
\put(120,105){\line(0,1){10}}
\put(150,110){\line(0,1){5}}
\put(15,92.5){\makebox(0,0){$\cdots$}}
\put(15,97.5){\makebox(0,0){$\cdots$}}
\put(15,102.5){\makebox(0,0){$\cdots$}}
\put(15,107.5){\makebox(0,0){$\cdots$}}
\put(15,112.5){\makebox(0,0){$\cdots$}}
\put(45,97.5){\makebox(0,0){$\cdots$}}
\put(45,102.5){\makebox(0,0){$\cdots$}}
\put(45,107.5){\makebox(0,0){$\cdots$}}
\put(45,112.5){\makebox(0,0){$\cdots$}}
\put(75,102.5){\makebox(0,0){$\cdots$}}
\put(75,107.5){\makebox(0,0){$\cdots$}}
\put(75,112.5){\makebox(0,0){$\cdots$}}
\put(105,107.5){\makebox(0,0){$\cdots$}}
\put(105,112.5){\makebox(0,0){$\cdots$}}
\put(135,112.5){\makebox(0,0){$\cdots$}}
\put(5,90){\line(0,1){25}}
\put(25,90){\line(0,1){25}}
\put(35,95){\line(0,1){20}}
\put(55,95){\line(0,1){20}}
\put(65,100){\line(0,1){15}}
\put(85,100){\line(0,1){15}}
\put(95,105){\line(0,1){10}}
\put(115,105){\line(0,1){10}}
\put(125,110){\line(0,1){5}}
\put(145,110){\line(0,1){5}}
\put(40,85){\vector(-1,0){40}}
\put(50,85){\vector(1,0){40}}
\put(45,85){\makebox(0,0){$\scriptstyle b$}}
\put(155,100){\makebox(0,0){$\scriptstyle (w)$}}
\end{picture}
\raisebox{15pt}{$\stackrel{\nabla}{\verylongrightarrow}\quad$}
\begin{picture}(170,35)(0,85)
\put(0,115){\line(1,0){150}}
\put(0,110){\line(1,0){150}}
\put(0,105){\line(1,0){120}}
\put(0,100){\line(1,0){90}}
\put(0,95){\line(1,0){60}}
\put(0,90){\line(1,0){30}}
\put(0,90){\line(0,1){25}}
\put(30,90){\line(0,1){25}}
\put(60,95){\line(0,1){20}}
\put(90,100){\line(0,1){15}}
\put(120,105){\line(0,1){10}}
\put(150,110){\line(0,1){5}}
\put(15,92.5){\makebox(0,0){$\cdots$}}
\put(15,97.5){\makebox(0,0){$\cdots$}}
\put(15,102.5){\makebox(0,0){$\cdots$}}
\put(15,107.5){\makebox(0,0){$\cdots$}}
\put(15,112.5){\makebox(0,0){$\cdots$}}
\put(45,97.5){\makebox(0,0){$\cdots$}}
\put(45,102.5){\makebox(0,0){$\cdots$}}
\put(45,107.5){\makebox(0,0){$\cdots$}}
\put(45,112.5){\makebox(0,0){$\cdots$}}
\put(75,102.5){\makebox(0,0){$\cdots$}}
\put(75,107.5){\makebox(0,0){$\cdots$}}
\put(75,112.5){\makebox(0,0){$\cdots$}}
\put(105,107.5){\makebox(0,0){$\cdots$}}
\put(105,112.5){\makebox(0,0){$\cdots$}}
\put(135,112.5){\makebox(0,0){$\cdots$}}
\put(5,90){\line(0,1){25}}
\put(25,90){\line(0,1){25}}
\put(35,95){\line(0,1){20}}
\put(55,95){\line(0,1){20}}
\put(65,100){\line(0,1){15}}
\put(85,100){\line(0,1){15}}
\put(95,105){\line(0,1){10}}
\put(115,105){\line(0,1){10}}
\put(125,110){\line(0,1){5}}
\put(145,110){\line(0,1){5}}
\put(40,85){\vector(-1,0){40}}
\put(60,85){\vector(1,0){35}}
\put(50,85){\makebox(0,0){$\scriptstyle b+1$}}
\put(155,100){\makebox(0,0){$\scriptstyle (w)\makebox[0pt][l]{$\,,$}$}}
\thicklines
\put(90,100){\line(1,0){5}}
\put(90,105){\line(1,0){5}}
\put(90,100){\line(0,1){5}}
\put(95,100){\line(0,1){5}}
\end{picture}\end{equation}
then $w=v+b+1-r$. More generally, taking conventions from
\cite{eProjective}, if $\nabla$ and
$\hat\nabla$ are two torsion-free connections in the same projective class
\[\hat\nabla_a\phi_b=\nabla_a\phi_b-\Upsilon_a\phi_b-\Upsilon_b\phi_a\]
for some $1$-form~$\Upsilon_a$, then for $\phi_{bc\cdots d}$ having symmetries
and projective weight specif\/ied by the left hand side of~(\ref{firstorder}),
\[\pi(\hat\nabla_a\phi_{bc\cdots d})=
\pi\big(\nabla_a\phi_{bc\cdots d}
+(w-(v+b+1-r))\Upsilon_a\phi_{bc\cdots d}\big),\]
where $\pi$ is the Young projector corresponding to the right hand side
of~(\ref{firstorder}). We may iterate this formula, adding more boxes to the
$r^{\mathrm{th}}$ row (assuming that there is room to do so) and, each time,
both $b$ and $v$ increase by one. Suppressing indices, after $s$
iterations the result is that
\[\pi(\hat\nabla^s\phi)=\pi\big(
(\nabla+(k-2s+2)\Upsilon\phi)\cdots
(\nabla+(k-4)\Upsilon\phi)(\nabla+(k-2)\Upsilon\phi)
(\nabla+k\Upsilon)\phi\big),\]
where we are writing $k=w-(v+b+1-r)$. In particular, if $s=w+r-v-b$ as it
is in the case of a BGG operator, then this iteration reads
\[\pi(\hat\nabla^s\phi)=\pi\big(
(\nabla-k\Upsilon\phi)\cdots
(\nabla+(k-4)\Upsilon\phi)(\nabla+(k-2)\Upsilon\phi)
(\nabla+k\Upsilon)\phi\big),\]
where $k=s-1$.
So far, this conclusion holds under any projective change of connection but now
we specialise to the case of the round connection $\nabla$ on the sphere, being
projectively equivalent to the f\/lat connection $\hat\nabla$ (under gnomonic
projection). The general formula \cite[(3.4)]{eProjective} for the
change in the Ricci tensor specialises to
\[0=g_{ab}-\nabla_a\Upsilon_b+\Upsilon_a\Upsilon_b\qquad\mbox{or, suppressing
indices,}\qquad\nabla\Upsilon=g+\Upsilon^2.\]
In this equation $\Upsilon$ is viewed as a tensor but if it is viewed as an
operator $\phi\stackrel{\Upsilon}{\longmapsto}\Upsilon\phi$, then we should
write $\nabla\Upsilon=\Upsilon\nabla+g+\Upsilon^2$. This equation allows us to
deal with the iterated formula above. For example, when $k=2$, also bearing in
mind that as operators $\nabla$ and $g$ commute,
 \begin{gather*}
 \underline{\pi\big(\hat\nabla^3\phi\big)}=
 \pi\big((\nabla-2\Upsilon)\nabla(\nabla+2\Upsilon)\big)\\
\hphantom{\underline{\pi\big(\hat\nabla^3\phi\big)}}{} = \pi\big(\nabla^3-2\Upsilon\nabla^2+2\nabla(\nabla\Upsilon)
-4\Upsilon(\nabla\Upsilon)\big)\\
\hphantom{\underline{\pi\big(\hat\nabla^3\phi\big)}}{} =
\pi\big(\nabla^3-2\Upsilon\nabla^2+2\nabla\big(\Upsilon\nabla+g+\Upsilon^2\big)
-4\Upsilon\big(\Upsilon\nabla+g+\Upsilon^2\big)\big)\\
\hphantom{\underline{\pi\big(\hat\nabla^3\phi\big)}}{} = \pi\big(\nabla^3-2\Upsilon\nabla^2+2(\nabla\Upsilon)\nabla+2g\nabla
+2(\nabla\Upsilon)\Upsilon
-4\Upsilon^2\nabla-4g\Upsilon-4\Upsilon^3\big)\\
\hphantom{\underline{\pi\big(\hat\nabla^3\phi\big)}}{} = \pi\big(\nabla^3-2\Upsilon\nabla^2+2\big(\Upsilon\nabla+g+\Upsilon^2\big)\nabla+2g\nabla
+2(\nabla\Upsilon)\Upsilon
-4\Upsilon^2\nabla-4g\Upsilon-4\Upsilon^3\big)\\
\hphantom{\underline{\pi\big(\hat\nabla^3\phi\big)}}{} = \pi\big(\nabla^3+4g\nabla+2\big(\Upsilon\nabla+g+\Upsilon^2\big)\Upsilon
-2\Upsilon^2\nabla-4g\Upsilon-4\Upsilon^3\big)\\
\hphantom{\underline{\pi\big(\hat\nabla^3\phi\big)}}{} = \pi\big(\nabla^3+4g\nabla+2\Upsilon(\nabla\Upsilon)
-2\Upsilon^2\nabla-2g\Upsilon-2\Upsilon^3\big)\\
\hphantom{\underline{\pi\big(\hat\nabla^3\phi\big)}}{} = \pi\big(\nabla^3+4g\nabla+2\Upsilon\big(\Upsilon\nabla+g+\Upsilon^2\big)
-2\Upsilon^2\nabla-2g\Upsilon-2\Upsilon^3\big)\\
\hphantom{\underline{\pi\big(\hat\nabla^3\phi\big)}}{} = \pi(\nabla^3+4g\nabla) =
\underline{\pi\big(\big(\nabla^2+4g\big)\nabla\big)},
\end{gather*}
as advertised in the statement of the theorem. For higher~$k$, direct
calculations rapidly get out of hand. Instead, it suf\/f\/ices to prove the
following lemma in which we have isolated the required algebra (and then we
prove the lemma by indirect means).
\end{proof}

\begin{remark}Our curvature normalisation (\ref{cc}) implies that the Ricci
tensor on our round sphere is given by
\[R_{ab}\equiv R_{ca}{}^c{}_b=(n-1)g_{ab}.\]
Thus, the metric $g_{ab}$ coincides with $\frac{1}{n-1}R_{ab}$, a tensor
generally known in projective dif\/ferential geometry \cite{eProjective} as the
{\em Schouten tensor\/} or {\em Rho-tensor\/}~$\Rho_{ab}$. Replacing $g$ by
$\Rho$ in the formul{\ae} of Theorem~\ref{round_formulae} gives expressions
that are valid on any space of constant curvature. More generally, there is a
Rho-tensor that arises in similar contexts~\cite{cds} within parabolic
dif\/ferential geometry~\cite{thebook}.
\end{remark}
\begin{lemma}\label{pretricky} Define an associative algebra
${\mathcal{R}}={\mathbb{R}}\langle\nabla,\Upsilon\rangle$ with generators
subject to the `Riccati relation'
$\nabla\Upsilon=\Upsilon\nabla+1+\Upsilon^2$. Then the following identities
hold in~${\mathcal{R}}$. If $k$ is even, then
\begin{gather*}
(\nabla-k\Upsilon)(\nabla-(k-2)\Upsilon)\cdots
(\nabla+(k-2)\Upsilon)(\nabla+k\Upsilon)=
(\nabla^2+k^2)\cdots(\nabla^2+16)(\nabla^2+4)\nabla\!
\end{gather*}
and, if $k$ is odd, then
\begin{gather*}
(\nabla-k\Upsilon)(\nabla-(k-2)\Upsilon)\cdots
(\nabla+(k-2)\Upsilon)(\nabla+k\Upsilon)=
(\nabla^2+k^2)\cdots(\nabla^2+9)(\nabla^2+1).
\end{gather*}
\end{lemma}
\begin{proof} The algebra ${\mathcal{R}}$ may be realised by the
following dif\/ferential operators on the circle
\[f(\theta)\stackrel{\nabla}{\longmapsto} df(\theta)/d\theta,\qquad
f(\theta)\stackrel{\Upsilon}{\longmapsto}(\tan\theta)f(\theta).\]
To see this, note that these operators certainly satisfy the Riccati
relation and we are required, therefore, to show that they satisfy no further
relations. Within ${\mathcal{R}}$ we may normalise any element as follows. By
induction, the Riccati relation extends to
\[\nabla\Upsilon^\ell=
\Upsilon^\ell\nabla+\ell\big(\Upsilon^{\ell-1}+\Upsilon^{\ell+1}\big),\qquad
\forall\,\ell\geq 1\]
whence
\[\nabla^k\Upsilon^\ell
=\nabla^{k-1}\big(\Upsilon^\ell\nabla
+\ell\big(\Upsilon^{\ell-1}+\Upsilon^{\ell+1}\big)\big)
=\big(\nabla^{k-1}\Upsilon^\ell\big)\nabla
+\ell\big(\nabla^{k-1}\Upsilon^{\ell-1}\big)+\ell\big(\nabla^{k-1}\Upsilon^{\ell+1}\big)\]
and it follows by induction on $k$ that
\begin{gather*}
\nabla^k\Upsilon^\ell=
\Upsilon^\ell\nabla^k+k\ell\big(\Upsilon^{\ell-1}+\Upsilon^{\ell+1}\big)\nabla^{k-1}\\
\hphantom{\nabla^k\Upsilon^\ell=}{}
+\tfrac12k(k-1)\ell\big((\ell-1)\Upsilon^{\ell-2}+2\ell\Upsilon^\ell
+(\ell+1)\Upsilon^{\ell+2}\big)\nabla^{k-2}+\cdots,
\end{gather*}
where the ellipsis $\cdots$ denotes terms of lower order in $\nabla$ with
coef\/f\/icients that are real polynomial in~$\Upsilon$. It follows that every
element of ${\mathcal{R}}$ can be written uniquely in the form
\[ \sum_{p=0}^{k}A_p(\Upsilon)\,\nabla^p\]
for suitable real polynomials $A_p(\Upsilon)$. In our claimed realisation,
such an element is represented by the dif\/ferential operator
\[ \sum_{p=0}^{k}A_p(\tan\theta)\,d^p/d\theta^p\]
and now it suf\/f\/ices to observe (by acting on
$1,\theta,\theta^2,\ldots,\theta^{k}$ near~$0$) that such a dif\/ferential
operator vanishes if and only if all the polynomials $A_p(\Upsilon)$ are zero.
(More precisely, we should restrict the action of such operators to smooth
functions on $(-\pi/2,\pi/2)$ or some other suitable function space.)

Having realised ${\mathcal{R}}$ by dif\/ferential operators, we are reduced to
proving identities amongst these operators. This is accomplished in the
following lemma.\end{proof}

\begin{lemma}\label{tricky} Let $D$ denote the differential operator
$f(\theta)\mapsto (\cos^2\theta)(df(\theta)/d\theta)$ on the circle. Then the
following identities hold. If $k$ is even, then
\[\frac{1}{\cos^k\theta}\frac{d}{d\theta}
\left(D^k\left(\frac{f(\theta)}{\cos^k\theta}\right)\right)
=\left(\frac{d^2}{d\theta^2}+k^2\right)\cdots
\left(\frac{d^2}{d\theta^2}+16\right)
\left(\frac{d^2}{d\theta^2}+4\right)\frac{d}{d\theta}f(\theta)\] and, if $k$ is
odd, then \[\frac{1}{\cos^k\theta}\frac{d}{d\theta}
\left(D^k\left(\frac{f(\theta)}{\cos^k\theta}\right)\right)
=\left(\frac{d^2}{d\theta^2}+k^2\right)\cdots
\left(\frac{d^2}{d\theta^2}+9\right)
\left(\frac{d^2}{d\theta^2}+1\right)f(\theta).\]
\end{lemma}
\begin{proof} Writing $D_{k+1}$ for the operators on the left-hand-sides of the
displays in this lemma,  the following identity is easily verif\/ied
\[
D_{k+3} = \left( \frac{d} { d\theta}
-(k + 2) \tan \theta\right)D_{k+1}\left( \frac{d}{ d \theta} + (k + 2) \tan \theta\right).\]
It follows that, in our realisation of the algebra ${\mathcal{R}}$,
we obtain the expressions on the left-hand-sides of the displays in
the Lemma~\ref{pretricky}.

On the other hand, the operators on the right-hand-sides of the
claimed identities, in the current lemma, are characterised up to scale
as annihilating the functions
\[\cos(k\theta),\quad \sin(k\theta),\quad \cos((k-2)\theta),\quad
\sin((k-2)\theta),\quad \cos((k-4)\theta),\quad \dots.\]
Since all operators have $d^{k+1}/d\theta^{k+1}$ as leading term, it is
therefore suf\/f\/icient to show that the left hand sides of these purported
identities have the same property. Notice that there is an invertible
relationship
\[\cos(m\theta)=2^{m-1}\cos^m\theta+\cdots=T_m(\cos\theta),\]
where $T_m$ is the $m^{\mathrm{th}}$ Chebyshev polynomial of the f\/irst kind
and a similar invertible relationship
\[\sin(m\theta)=(\sin\theta)\big(2^{m-2}\cos^{m-1}\theta+\cdots\big)
=(\sin\theta)U_{m-1}(\cos\theta)\]
where $U_{m-1}$ is the $(m-1)^{\mathrm{st}}$ Chebyshev polynomial of the
second kind. Only the degree of these Chebyshev polynomials concerns us and
it now suf\/f\/ices to show that $D^{k+1}$ annihilates the f\/irst $k+1$ of
\[1,\quad\frac{\sin\theta}{\cos\theta},\quad
\frac{1}{\cos^2\theta},\quad\frac{\sin\theta}{\cos^3\theta},\quad
\frac{1}{\cos^4\theta},\quad\frac{\sin\theta}{\cos^5\theta},\quad
\frac{1}{\cos^6\theta},\quad \frac{\sin\theta}{\cos^7\theta},\quad\dots.\]
Since
\[D\frac{1}{\cos^{2\ell}\theta}=2\ell\frac{\sin\theta}{\cos^{2\ell-1}\theta}
\qquad\mbox{and}\qquad
D\frac{\sin\theta}{\cos^{2\ell-1}\theta}
=(2\ell-1)\frac{1}{\cos^{2(\ell-1)}\theta}
-2(\ell-1)\frac{1}{\cos^{2(\ell-2)}\theta},\]
this follows easily by induction.
\end{proof}
\begin{remark} More generally, the formula
\[D_{k+1}f=u^{-k-2}\left(u^2\frac{d}{dx}\right)^{k+1}u^{-k}f,
\qquad\mbox{where} \quad \left(\frac{d^2}{dx^2}+\Phi\right)u=0\]
is used in~\cite{g} to derive expressions for BGG operators in conformal
geometry and in~\cite{cds} these expressions are extended to parabolic
geometries in general. Lemma~\ref{tricky} concerns the case $\Phi=1$.
\end{remark}

\section{An example}

The following is amongst the simplest of non-trivial examples. Let us take
$n=2$ so $G={\mathrm{SL}}(3,{\mathbb{R}})$ and let
${\mathbb{V}}={\mathbb{R}}^3$, regarded as column vectors with
${\mathrm{SL}}(3,{\mathbb{R}})$ acting by left matrix multiplication. We saw in~(\ref{Hstandard}) how $H$ splits ${\mathbb{R}}^3$ and in Section~\ref{filtering} that
this results in the f\/iltering
\[{\mathbb{R}}^3={\mathbb{R}}^2+{\mathbb{R}}\quad\mbox{as a $P$-module}\]
where ${\mathrm{SL}}(2,{\mathbb{R}})\subset P$ acts on ${\mathbb{R}}^2$ as the
standard representation and acts trivially on~${\mathbb{R}}$. Dropping
projective weights, the corresponding bundle $V$ on ${\mathbb{RP}}_2$ is
f\/iltered
\[V=V_0+V_1=T+\Lambda^0,\]
where $T$ is the tangent bundle and $\Lambda^0$ is the trivial bundle. The
$E_0$-level (\ref{E0}) of our spectral sequence becomes
\[\begin{picture}(300,70)(0,15)
\put(0,20){\vector(1,0){190}}
\put(0,20){\vector(0,1){60}}
\put(183,14){\makebox(0,0){$p$}}
\put(-4,71){\makebox(0,0){$q$}}
\put(10,60){\makebox(0,0){$T$}}
\put(50,60){\makebox(0,0){$\Lambda^1\otimes T$}}
\put(100,60){\makebox(0,0){$\Lambda^2\otimes T$}}
\put(150,60){\makebox(0,0){$0$}}
\put(50,45){\makebox(0,0){$\uparrow\scriptstyle\partial$}}
\put(100,45){\makebox(0,0){$\uparrow\scriptstyle\partial$}}
\put(10,30){\makebox(0,0){$0$}}
\put(50,30){\makebox(0,0){$\Lambda^0$}}
\put(100,30){\makebox(0,0){$\Lambda^1$}}
\put(150,30){\makebox(0,0){$\Lambda^2$}}
\end{picture}\]
and one checks that $\partial:\Lambda^0\to\Lambda^1\otimes T$ and
$\partial:\Lambda^1\to\Lambda^2\otimes T$ are given by
\[\mu\mapsto\delta_b{}^c\mu\qquad\mbox{and}\qquad
\mu_b\mapsto\delta_{[a}{}^c\mu_{b]},\]
respectively. In this simple case, we do not need Kostant's Theorem to see that
 \begin{gather*}
\partial: \ \Lambda^0\to\Lambda^1\otimes T\enskip\mbox{is injective with
cokernel}=(\Lambda^1\otimes T)_\circ,\\
\partial: \ \Lambda^1\to\Lambda^2\otimes T\enskip\mbox{is an isomorphism},
\end{gather*}
where $(\Lambda^1\otimes T)_\circ$ denotes the trace-free part of
$\Lambda^1\otimes T$. Therefore, the $E_1$-level (\ref{sparse}) is
\[\begin{picture}(300,70)(0,15)
\put(0,20){\vector(1,0){190}}
\put(0,20){\vector(0,1){60}}
\put(183,14){\makebox(0,0){$p$}}
\put(-4,71){\makebox(0,0){$q$}}
\put(10,60){\makebox(0,0){$T$}}
\put(25,60){\makebox(0,0){$\to$}}
\put(60,60){\makebox(0,0){$(\Lambda^1\otimes T)_\circ$}}
\put(110,60){\makebox(0,0){$0$}}
\put(150,60){\makebox(0,0){$0$}}
\put(10,30){\makebox(0,0){$0$}}
\put(60,30){\makebox(0,0){$0$}}
\put(110,30){\makebox(0,0){$0$}}
\put(150,30){\makebox(0,0){$\Lambda^2$}}
\end{picture}\]
and we obtain
\[0\to{\mathbb{R}}^3
\to T\to (\Lambda^1\otimes T)_\circ\to \Lambda^2\to 0\]
as the resulting BGG complex. Writing out the f\/lat connection $\nabla$ on
$V$ in terms of the round connection gives
\[\left\lgroup\begin{array}c\sigma^c\\ \mu\end{array}\right\rgroup
\stackrel{\nabla}{\longmapsto}
\left\lgroup\begin{array}c\nabla_b\sigma^c+\delta_b{}^c\mu\\
\nabla_b\mu-\sigma_b\end{array}\right\rgroup\qquad\mbox{and}\qquad
\left\lgroup\begin{array}c\sigma_b{}^c\\ \mu_b\end{array}\right\rgroup
\stackrel{\nabla}{\longmapsto}
\left\lgroup\begin{array}c\nabla_{[a}\sigma_{b]}{}^c+\delta_{[a}{}^c\mu_{b]}\\
\nabla_{[a}\mu_{b]}+\sigma_{[ab]}\end{array}\right\rgroup\]
for the two operators in $V\to\Lambda^1\otimes V\to\Lambda^2\otimes V$ and,
noting that
\[\Lambda^1\otimes V\ni\left\lgroup\begin{array}c\sigma_b{}^c\\
-\nabla_c\sigma_b{}^c\end{array}\right\rgroup
\stackrel{\nabla}{\longmapsto}
-\left\lgroup\begin{array}c0\\
\nabla_{c}\nabla_{[a}\sigma_{b]}{}^c+g_{c[a}\sigma_{b]}{}^c
\end{array}\right\rgroup,\]
we obtain a formula for the projectively invariant operator
$(\Lambda^1\otimes T)_\circ\stackrel{\nabla^{(2)}}{\longrightarrow}\Lambda^2$
in agreement with Theorem~\ref{round_formulae}.

\renewcommand{\thesection}{Appendix:}
\section[Lie algebra cohomology as a geometrical construction]{{}\hspace*{-5pt}Lie algebra cohomology as a geometrical construction}

Although in this article we shall need only the cohomology of an Abelian Lie
algebra, we take the opportunity here to describe the cohomology of a general
Lie algebra ${\mathfrak{g}}$ in terms of dif\/ferential geometry on~$G$, a Lie
group whose Lie algebra is~${\mathfrak{g}}$. We believe that for a general
parabolic geometry, we shall need this geometric interpretation for a nilpotent
Lie algebra. Suppose ${\mathbb{V}}$ is a $G$-module and use the same notation
for the corresponding representation of~${\mathfrak{g}}$. Following but
adapting~\cite{ce}, we are going to present the Lie algebra cohomology
$H^r({\mathfrak{g}},{\mathbb{V}})$ as a geometrical construction on~$G$. Beware
that $G$ is no longer the Lie group ${\mathrm{SL}}(n+1,{\mathbb{R}})$ as it was
until now. This section is written to be self-contained with the aim of being
useful elsewhere. This material is well-known to experts and implicit
in~\cite{ce} but we believe it worthwhile laying out the details.

We shall view $G$ as a homogeneous space under its own action on the left. Its
tangent bundle~$TG$ is then regarded as a homogeneous bundle and can be
identif\/ied as $G\times{\mathfrak{g}}$, where ${\mathfrak{g}}$ is the Lie
algebra of $G$. It is convenient to write this isomorphism as
\begin{equation}\label{m-c}
X\mapsto X\intprod\theta \qquad\mbox{for vector f\/ields $X$ on $G$},
\end{equation}
where $\theta$ is a $1$-form on $G$ with values in~${\mathfrak{g}}$ known as
the {\em Maurer--Cartan\/} form~\cite{ss}. To compute with $\theta$ it is
convenient to write functions on $G$ with values in ${\mathfrak{g}}$ as
$X^\alpha$ and then (\ref{m-c}) becomes
\[
X^a\mapsto X^\alpha\equiv\theta_a^\alpha X^a\qquad \mbox{with inverse}
\quad X^\alpha\mapsto X^a\equiv\phi_\alpha^a X^\alpha,
\]
where $\phi_\alpha^a\theta^\beta_a=\delta_\alpha{}^\beta$ and $\theta_a^\alpha\phi_\alpha^b=\delta_a{}^b$.
A vector f\/ield $X^a$ on $G$ is {\em left-invariant\/} if and only if the
corresponding function $X^\alpha:G\to{\mathfrak{g}}$ is constant. Choosing any
torsion-free af\/f\/ine connection $\nabla_a$ on $G$ and expanding the def\/inition
$\nabla_a(\theta_b^\alpha X^b)=0$, the left-invariant vector f\/ields are those
that satisfy
\[\theta_b^\alpha \nabla_aX^b-X^b \nabla_a\theta_b^\alpha=0\]
from which it follows easily that the Lie bracket of two left-invariant
vector f\/ields is again left-invariant.
Since the left-invariant vector f\/ields on $G$ are of the form $\phi(X)$ for
$X\in{\mathfrak{g}}$ we may def\/ine the Lie bracket
on ${\mathfrak{g}}$ by transportation: 
\begin{equation}\label{defnLie}
\phi([X,Y])=[\phi(X),\phi(Y)]\qquad\mbox{for} \quad X,Y\in{\mathfrak{g}}.
\end{equation}
For computational purposes, let us write
$[X,Y]^\gamma=\Gamma_{\alpha\beta}{}^\gamma X^\alpha Y^\beta$ for the Lie
bracket on~${\mathfrak{g}}$. Then we can write out~(\ref{defnLie}) explicitly
as
 \begin{gather*}
\phi_\gamma^c\Gamma_{\alpha\beta}{}^\gamma X^\alpha Y^\beta =
[\phi_\alpha^a X^\alpha,\phi_\beta^b Y^\beta]^c
 = \phi_\alpha^a X^\alpha \nabla_a(\phi_\beta^c Y^\beta)-
\phi_\alpha^b Y^\alpha \nabla_b(\phi_\beta^c X^\beta)\\
\hphantom{\phi_\gamma^c\Gamma_{\alpha\beta}{}^\gamma X^\alpha Y^\beta =
[\phi_\alpha^a X^\alpha,\phi_\beta^b Y^\beta]^c}{}
 = (\phi_\alpha^a \nabla_a\phi_\beta^c-\phi_\beta^a \nabla_a\phi_\alpha^c)
X^\alpha Y^\beta
\end{gather*}
or, in other words, as
\[
\Gamma_{\alpha\beta}{}^\gamma\phi_\gamma^c
=\phi_\alpha^a \nabla_a\phi_\beta^c-\phi_\beta^a \nabla_a\phi_\alpha^c.
\]
Bearing in mind that $\phi_\beta^b\theta_b^\gamma=\delta_\beta{}^\gamma$ whence
$\theta_b^\gamma \nabla_a\phi_\beta^b+\phi_\beta^b \nabla_a\theta_b^\gamma=0$,
we may rewrite this as
\[\Gamma_{\alpha\beta}{}^\gamma
=-\phi_\alpha^a\phi_\beta^b \nabla_a\theta_b^\gamma+
\phi_\beta^a\phi_\alpha^b \nabla_a\theta_b^\gamma
=-\phi_\alpha^a\phi_\beta^b(\nabla_a\theta_b^\gamma-\nabla_b\theta_a^\gamma)\]
and f\/inally as
\begin{equation}\label{m-c-eq1}
\nabla_a\theta_b^\gamma-\nabla_b\theta_a^\gamma
+\Gamma_{\alpha\beta}{}^\gamma\theta_a^\alpha\theta_b^\beta=0.
\end{equation}
This formula employs an arbitrary torsion-free connection. Without indices and
without this connection  it is more usually written as
\begin{equation}\label{m-c-eq2}  d\theta+\tfrac12[\theta,\theta]=0\qquad
\mbox{or}\qquad d\theta+\theta\wedge\theta=0.
\end{equation}
In other words, the def\/inition (\ref{defnLie}) of the Lie bracket on
${\mathfrak{g}}$ is equivalent to (\ref{m-c-eq1}) or~(\ref{m-c-eq2}), usually
known as the {\em Maurer--Cartan\/} equation~\cite{ss}.

Now suppose ${\mathbb{V}}$ is a $G$-module and use the same notation for the
corresponding ${\mathfrak{g}}$-module. There are two canonically def\/ined
connections on the vector bundle $V=G\times{\mathbb{V}}$ over~$G$. One of them
is the evident f\/lat connection, ignoring the action of $G$ on ${\mathbb{V}}$.
We shall denote it by $d$ since it is the exterior derivative acting on
functions with values in~${\mathbb{V}}$. The other one takes the isomorphism
from (\ref{twist})
\begin{equation}\label{twistagain}
V=G\times{\mathbb{V}}\cong G\times{\mathbb{V}}\qquad\mbox{by}\quad
(g,v)\mapsto(g,gv)\end{equation}
and pulls back the evident f\/lat connection on the right hand side as was done
in Section~\ref{outline} and we shall denote this one by $\nabla$ as was done there.
To relate these two connections more explicitly suppose $f:G\to{\mathbb{V}}$ is
a section of $V$ that is constant after the twisting~(\ref{twistagain}). It
means that the function $g\mapsto gf(g)$ is constant. If so, then for f\/ixed
$g\in G$ and~$X\in{\mathfrak{g}}$, the function
\[{\mathbb{R}}\ni t\mapsto ge^{tX}f\big(ge^{tX}\big)\in{\mathbb{V}}\]
is constant. Equivalently, the function $t\mapsto e^{tX}f(ge^{tX})$ is
constant and so
\[0=\frac{d}{dt}\big(e^{tX}f\big(ge^{tX}\big)\big)\Big|_{t=0}
=\frac{d}{dt}f\big(ge^{tX}\big)\big|_{t=0}+Xf(g)=(\phi(X)f+Xf)(g),\]
where this last equality is due to the f\/low of the left-invariant vector f\/ield
$\phi(X)$ being the one-parameter subgroup of right-translations
$g\mapsto ge^{tX}$ (see, e.g.~\cite{w}). For computational purposes, let us
write $X^\alpha\mapsto X^\alpha\rho_\alpha$ where
$\rho_\alpha\in{\mathfrak{g}}^*\otimes\End({\mathbb{V}})$ for the action of
${\mathfrak{g}}$ on~${\mathbb{V}}$. Then $g\mapsto gf(g)$ is constant if and
only if
\[0=\phi(X)f+Xf=X^\alpha\phi_\alpha^ad_af+X^\alpha\rho_\alpha f
=X^a(d_af+\theta_a^\alpha\rho_\alpha f)\]
for all left-invariant vector f\/ields~$X^a$. It follows that the connection
$\nabla_a$ on $V$ is given by
\begin{equation}\label{nabla}
f\mapsto\nabla_af=d_af+\theta_a^\alpha\rho_\alpha f\qquad
\mbox{or, without indices, as}\quad
f\mapsto\nabla f=df+\theta f,\end{equation}
where $\theta\in\Lambda^1\otimes{\mathfrak{g}}$ is the Maurer--Cartan form. As a
check, the dif\/ferential in the coupled de~Rham sequence (\ref{coupled}) is
\[\Lambda^p\otimes{\mathbb{V}}\ni\omega\stackrel{\nabla}{\longmapsto}
d\omega+\theta\wedge\omega
\in\Lambda^{p+1}\otimes{\mathbb{V}}\]
and the Maurer--Cartan equation~(\ref{m-c-eq2}) shows that the composition
${\mathbb{V}}\stackrel{\nabla}{\to}\Lambda^1\otimes{\mathbb{V}}
\stackrel{\nabla}{\to}\Lambda^2\otimes{\mathbb{V}}$ is given by
\[\nabla^2f=d(df+\theta f)+\theta\wedge(df+\theta f)
=(d\theta+\theta\wedge\theta)f=0\]
and the connection $\nabla$ is f\/lat, as expected. More generally, the whole
sequence
\begin{equation}\label{coupledtobbV}
0\to{\mathbb{V}}\stackrel{\nabla}{\longrightarrow}
\Lambda^1\otimes{\mathbb{V}}\stackrel{\nabla}{\longrightarrow}
\Lambda^2\otimes{\mathbb{V}}
\stackrel{\nabla}{\longrightarrow}\cdots\stackrel{\nabla}{\longrightarrow}
\Lambda^p\otimes{\mathbb{V}}
\stackrel{\nabla}{\longrightarrow}
\Lambda^{p+1}\otimes{\mathbb{V}}
\stackrel{\nabla}{\longrightarrow}
\cdots\end{equation}
is a complex.

Consider the linear mapping
\[\Lambda^p{\mathfrak{g}}^*\otimes{\mathbb{V}}\ni v_{\alpha\beta\cdots\gamma}
\stackrel{\theta}{\longmapsto}v_{ab\cdots c}\equiv
\theta_a^\alpha\theta_b^\beta\cdots\theta_c^\gamma
v_{\alpha\beta\cdots\gamma}\in\Gamma(G,\Lambda^p\otimes{\mathbb{V}})\]
in which ${\mathbb{V}}$ is just a passenger (i.e.~plays no r\^{o}le). We shall
refer to the resulting ${\mathbb{V}}$-valued $p$-form as {\em left-invariant\/}
just as we would if ${\mathbb{V}}$ were absent.
\begin{lemma} The connection
$\nabla:\Lambda^p\otimes{\mathbb{V}}\to\Lambda^{p+1}\otimes{\mathbb{V}}$
preserves left-invariance.\end{lemma}
\begin{proof}We use the Maurer--Cartan equation (\ref{m-c-eq1}) to compute
\begin{gather}
\nabla_{[a}v_{bc\cdots d]} =d_{[a}(\theta_b^\beta\theta_c^\gamma\cdots\theta_{d]}^\delta
v_{\beta\gamma\cdots\delta})+
\theta_{[a}^\alpha\theta_b^\beta\theta_c^\gamma\cdots\theta_{d]}^\delta
\rho_\alpha v_{\beta\gamma\cdots\delta}
\nonumber\\
 \hphantom{\nabla_{[a}v_{bc\cdots d]}}{} = p(d_{[a}\theta_b^\epsilon)\theta_c^\gamma\cdots\theta_{d]}^\delta
v_{\epsilon\gamma\cdots\delta}+
\theta_{[a}^\alpha\theta_b^\beta\theta_c^\gamma\cdots\theta_{d]}^\delta
\rho_\alpha v_{\beta\gamma\cdots\delta}
\nonumber\\
 \hphantom{\nabla_{[a}v_{bc\cdots d]}}{} =-(p/2)\Gamma_{\alpha\beta}{}^\epsilon
\theta_{[a}^\alpha\theta_b^\beta\theta_c^\gamma\cdots\theta_{d]}^\delta
v_{\epsilon\gamma\cdots\delta}+
\theta_{[a}^\alpha\theta_b^\beta\theta_c^\gamma\cdots\theta_{d]}^\delta
\rho_\alpha v_{\beta\gamma\cdots\delta}
\nonumber\\
 \hphantom{\nabla_{[a}v_{bc\cdots d]}}{} =\theta_{[a}^\alpha\theta_b^\beta\theta_c^\gamma\cdots\theta_{d]}^\delta
\rho_\alpha v_{\beta\gamma\cdots\delta}
+(-1)^p(p/2)
\theta_{[a}^\alpha\theta_b^\beta\theta_c^\gamma\cdots\theta_{d]}^\delta
\Gamma_{\alpha\beta}{}^\epsilon v_{\gamma\cdots\delta\epsilon}
\nonumber\\
 \hphantom{\nabla_{[a}v_{bc\cdots d]}}{} = \theta_a^\alpha\theta_b^\beta\theta_c^\gamma\cdots\theta_d^\delta
\left(\rho_{[\alpha}v_{\beta\gamma\cdots\delta]}
+(-1)^p(p/2)
\Gamma_{[\alpha\beta}{}^\epsilon v_{\gamma\cdots\delta]\epsilon}\right),\label{key}
\end{gather}
as required.
\end{proof}

In fact (\ref{key}) shows that $\nabla\theta v=\theta\partial v$, where
\begin{gather}
\partial: \ \Lambda^p{\mathfrak{g}}^*\otimes{\mathbb{V}}\to
\Lambda^{p+1}{\mathfrak{g}}^*\otimes{\mathbb{V}}
\qquad\mbox{is given by}\nonumber\\
\qquad
v_{\beta\gamma\cdots\delta}\mapsto\rho_{[\alpha}v_{\beta\gamma\cdots\delta]}
+(-1)^p(p/2)\Gamma_{[\alpha\beta}{}^\epsilon v_{\gamma\cdots\delta]\epsilon}.\label{thisispartial}
\end{gather}
It also follows that
\begin{equation}\label{Koszul}
0\to{\mathbb{V}}\stackrel{\partial}{\longrightarrow}
{\mathfrak{g}}^*\otimes{\mathbb{V}}\stackrel{\partial}{\longrightarrow}
\Lambda^2{\mathfrak{g}}^*\otimes{\mathbb{V}}
\stackrel{\partial}{\longrightarrow}\cdots\stackrel{\partial}{\longrightarrow}
\Lambda^p{\mathfrak{g}}^*\otimes{\mathbb{V}}
\stackrel{\partial}{\longrightarrow}
\Lambda^{p+1}{\mathfrak{g}}^*\otimes{\mathbb{V}}
\stackrel{\partial}{\longrightarrow}
\cdots\end{equation}
is a complex of~${\mathfrak{g}}$-modules. Alternatively, this may be directly
verif\/ied from (\ref{thisispartial}) using
\begin{itemize}\itemsep=0pt
\item
$\rho_{[\alpha}\rho_{\beta]}=\frac12\Gamma_{\alpha\beta}{}^\gamma\rho_\gamma$
(i.e.~that $\rho:{\mathfrak{g}}\to\End({\mathbb{V}})$ is a representation),
\item $\Gamma_{[\alpha\beta}{}^\delta\Gamma_{\gamma]\delta}{}^\epsilon=0$
(i.e.~the Jacobi identity in ${\mathfrak{g}}$),
\item $(Xv)_{\beta\gamma\cdots\delta}
=X^\alpha\rho_\alpha v_{\beta\gamma\cdots\delta}
+(-1)^ppX^\alpha\Gamma_{\alpha[\beta}{}^\epsilon
v_{\gamma\cdots\delta]\epsilon}$ (the action of ${\mathfrak{g}}$ on
$\Lambda^p{\mathfrak{g}}^*\otimes{\mathbb{V}}$).
\end{itemize}
We have shown that there are two ways of def\/ining Lie algebra cohomology as
follows.
\begin{theorem}\label{geometricrealisation} The Lie algebra cohomology
$H^r({\mathfrak{g}},{\mathbb{V}})$ may be defined as either
\begin{itemize}\itemsep=0pt
\item the cohomology of the Koszul complex~\eqref{Koszul}, or
\item the cohomology of the complex \eqref{coupledtobbV} restricted to
left-invariant forms.
\end{itemize}
\end{theorem}
\begin{remark} Although we use this theorem in the main body of this article,
it is easily avoided. In tackling a general parabolic geometry,
however, we believe that Theorem~\ref{geometricrealisation} will be essential.
\end{remark}

\begin{remark}
As a minor variation on this construction, suppose $G$ is enlarged to~$Q$, a
semi-direct product
\[Q=G_0\ltimes G\qquad\mbox{i.e.}\qquad
\begin{picture}(0,0)
\put(0,0){
${\mathrm{Id}}\to G\lhd Q\longrightarrow G_0\to{\mathrm{Id}}$}
\qbezier (67,-0.5) (74,-4) (81,-0.5)
\put(65,0){\vector(-3,1){0}}
\end{picture}\]
and suppose that ${\mathbb{V}}$ extends to a representation of~$Q$. We identify
$G$ with the $Q$-homogeneous space~$Q/G_0$, noting that when the action of $Q$
on $G=Q/G_0$ is restricted to~$G$ it coincides with its usual action of $G$ on
itself by left translation. The $Q$-homogeneous bundle
$V\equiv Q\times_{G_0}{\mathbb{V}}$ on~$Q/G_0$ is equipped with a f\/lat
connection $\nabla$ by dint of the canonical trivialisation
\[V=Q\times_{G_0}{\mathbb{V}}\ni(q,v)\mapsto(qG_0,qv)\in
Q/G_0\otimes{\mathbb{V}}=G\times{\mathbb{V}},\]
which clearly coincides with $\nabla$ def\/ined by~(\ref{twistagain}). This
connection is $Q$-equivariant. Consequently, not only does the twisted
de~Rham complex
\[V\xrightarrow{\,\nabla\,}\Lambda^1\otimes V
\xrightarrow{\,\nabla\,}\Lambda^2\otimes V
\xrightarrow{\,\nabla\,}\cdots
\xrightarrow{\,\nabla\,}\Lambda^p\otimes V
\xrightarrow{\,\nabla\,}\Lambda^{p+1}\otimes V
\xrightarrow{\,\nabla\,}\cdots\]
coincide with (\ref{coupledtobbV}) and thereby compute the Lie algebra
cohomology $H^r({\mathfrak{g}},{\mathbb{V}})$ when restricted to $G$-invariant
forms, but also the complex (\ref{Koszul}) is automatically one of
$Q$-modules where the $Q$-action on ${\mathbb{V}}$
is as supposed and the $Q$-action on ${\mathfrak{g}}^*$ is induced by the
conjugation action of $Q$ on $G$ (bearing in mind that $G$ is a normal
subgroup of~$Q$).
\end{remark}

\renewcommand{\thesubsection}{Addendum:}

\subsection[A canonical connection on $G$]{{}\hspace*{-3pt}A canonical connection on $\boldsymbol{G}$}

Again, although it is unnecessary for the current article and already known to
experts, we suspect that the following optional extra will be invaluable in
dealing with a general parabolic geometry. Since we already have established
suitable notation in the Appendix above, we take the opportunity of presenting
it here. Our canonical connection $D_a$ was introduced in~\cite{cs} as the
`$(+)$-connection' and Lemma~\ref{torsion_is_bracket} is stated without
proof as~\cite[Proposition~2.12]{knII}.

The trivialisation $T^*G=G\times{\mathfrak{g}}^*$ provided by the Maurer--Cartan
form also equips $G$ with a~canonical f\/lat af\/f\/ine connection $D_a$ def\/ined
by
\[D_a\omega_b\equiv\theta_b^\beta d_a\omega_\beta,\]
where $\omega_\beta\equiv\phi_\beta^c\omega_c$ and
$d_a$ on the right hand side of this equation simply takes the
gradient of a~function with values in~${\mathfrak{g}}^*$. If we expand using
any torsion-free af\/f\/ine connection $\nabla_a$
\[D_a\omega_b=\theta_b^\beta \nabla_a(\phi_\beta^c\omega_c)
=\nabla_a\omega_b+(\theta_b^\beta \nabla_a\phi_\beta^c)\omega_c
=\nabla_a\omega_b-(\phi_\beta^c\nabla_a\theta_b^\beta)\omega_c,\]
then we see that, for $f$ a smooth function,
\[D_aD_bf-D_bD_af
=(-\phi_\beta^c \nabla_a\theta_b^\beta+\phi_\beta^c \nabla_b\theta_a^\beta)D_cf
=-\phi_\gamma^c(\nabla_a\theta_b^\gamma-\nabla_b\theta_a^\gamma)D_cf\]
and so the canonical connection $D_a$ has torsion
\[T_{ab}{}^c=-(\nabla_a\theta_b^\gamma-\nabla_b\theta_a^\gamma)\phi_\gamma^c.\]
Alternatively, from (\ref{m-c-eq1}) we see that
\[T_{ab}{}^c\theta_c^\gamma=-(\nabla_a\theta_b^\gamma-\nabla_b\theta_a^\gamma)
=\Gamma_{\alpha\beta}{}^\gamma\theta_a^\alpha\theta_b^\beta,\]
which we record as the following lemma.
\begin{lemma}\label{torsion_is_bracket}
The torsion of\/ $D_a$ coincides with
the Lie bracket on\/~${\mathfrak{g}}$ under the Maurer--Cartan parallelism.
\end{lemma}
Notice that $D_a\theta_b^\beta=0$. It is another way to
characterise~$D_a$ and, indeed, is the main point of this construction as
follows.
\begin{lemma} Even locally, the kernel of the induced operator
\[D: \ \Lambda^p\to\Lambda^1\otimes\Lambda^p\]
is the left-invariant forms on~$G$.
\end{lemma}
\begin{proof}
Recall that the left-invariant forms are obtained as
\[\Lambda^p{\mathfrak{g}}^*\ni v_{\alpha\beta\cdots\gamma}
\stackrel{\theta}{\longmapsto}v_{ab\cdots c}\equiv
\theta_a^\alpha\theta_b^\beta\cdots\theta_c^\gamma
v_{\alpha\beta\cdots\gamma}\in\Gamma(G,\Lambda^p)\]
and it clear that such forms are annihilated by~$D$. Conversely, since $D$ is
f\/lat and all covariant constant sections are already accounted for, there can
be no more, even locally.
\end{proof}
\begin{remark} Of course, this lemma also holds for ${\mathbb{V}}$-valued
dif\/ferential forms where the connection is trivially coupled
with~${\mathbb{V}}$ and it is this that we have in mind in constructing
the BGG complex in general.
\end{remark}

\subsection*{Acknowledgements}
Eastwood is supported by the Australian Research Council.
 Gover is partly supported by the Royal Society of New Zealand
(Marsden Grant 10-UOA-113).
 The authors would like to thank Katharina Neusser for many valuable
conversations and the referees for their useful corrections and suggestions.

\pdfbookmark[1]{References}{ref}
 \LastPageEnding

\end{document}